\documentclass{amsart}

\newtheorem{theorem}{Theorem}[section]
\newtheorem{lemma}[theorem]{Lemma}

\newtheorem{proposition}[theorem]{Proposition}

\usepackage{graphicx}

\theoremstyle{definition}

\theoremstyle{remark}

\numberwithin{equation}{section}

\begin{document}

\title[Boundary structure of $3$-manifolds]
  {Boundary structure of hyperbolic $3$-manifolds admitting annular and toroidal fillings at large distance}

\author{Sangyop Lee}
\address{School of Mathematics, Korea Institute for Advanced Study, 207-42 Cheongryangridong, Dongdaemun-gu,
Seoul 130-012, Korea}
\email{slee@kias.re.kr}
%
%
\author{Masakazu Teragaito}
\address{Department of Mathematics and Mathematics Education, Hiroshima University,
1-1-1 Kagamiyama, Higashi-hiroshima, Japan 739-8524}
\email{teragai@hiroshima-u.ac.jp}
\thanks{
The second author was partially supported by Japan Society for the Promotion of Science,
Grant-in-Aid for Scientific Research (C), 16540071.
}%





\begin{abstract}
For a hyperbolic $3$-manifold $M$ with a torus boundary component,
all but finitely many Dehn fillings yield hyperbolic $3$-manifolds.
In this paper, we will focus on the situation where 
$M$ has two exceptional Dehn fillings: an annular filling and a toroidal filling.
For such situation, Gordon gave an upper bound $5$ for the distance between such slopes.
Furthermore, the distance $4$ is realized only by two specific manifolds, and $5$
is realized by a single manifold.
These manifolds all have a union of two tori as their boundaries.
Also, there is a manifold with three tori as its boundary which realizes the distance three.
We show that if the distance is three then the boundary of the manifold consists of at most three tori.
\end{abstract}

\maketitle

\section{Introduction}\label{sec:1}

Let $M$ be a hyperbolic $3$-manifold with a torus boundary component $T_0$.
For a slope $\gamma$ on $T_0$, $M(\gamma)$ denotes the manifold obtained by $\gamma$-Dehn filling on $M$.
That is, $M(\gamma)=M\cup V_\gamma$, where $V_\gamma$ is a solid torus,
glued to $M$ along $T_0$ in such a way that $\gamma$ bounds a disk in $V_\gamma$.
A $3$-manifold is said to be \textit{annular\/} (resp.\ \textit{toroidal\/}) 
if it contains an essential annulus (resp.\ torus).
Suppose that $M(\alpha)$ is annular and $M(\beta)$ is toroidal for slopes $\alpha$ and $\beta$ on $T_0$.
Gordon \cite{Go2} showed that $\Delta(\alpha,\beta)\le 5$, where
$\Delta(\alpha,\beta)$ denotes the \textit{distance\/} between two slopes, which is their minimal geometric intersection number.
Furthermore, Gordon and Wu \cite{GW} showed that the distance $5$ is realized by a single manifold and the distance $4$
is realized by two specific manifolds.
These manifolds are the exteriors of the Whitehead sister (or $(-2,3,8)$-pretzel) link, the Whitehead link and
the $2$-bridge link corresponding to $3/10$, using Conway's notation, in the $3$-sphere $S^3$.
Following Gordon \cite{Go4}, let us define
\begin{equation*}
\begin{split}
\Delta^k(A,T) &= \mbox{max}\{\Delta(\alpha,\beta)\ :\ \mbox{there is a hyperbolic $3$-manifold $M$ such that $\partial M$} \\
              & \quad\mbox{is a disjoint union of $k$ tori, and $\alpha,\beta$ are slopes on some component} \\
              & \quad\mbox{of $\partial M$ such that $M(\alpha)$ is annular and $M(\beta)$ is toroidal}\}\\
\end{split}
\end{equation*}
for $k\ge 2$.
($\Delta^k(X,Y)$ is defined similarly for other types $X,Y\in \{S,D,A,T\}$, but we do not need it.)
Thus $\Delta^2(A,T)=5$.
Also, there are infinitely many hyperbolic manifolds realizing the distance three \cite{GW}.
Among them, there is a hyperbolic $3$-manifold, called the \textit{magic manifold}, which is the exterior
of a certain $3$-component link in $S^3$.
Hence $\Delta^3(A,T)=3$.
Gordon \cite{Go4} gave an example showing $\Delta^k(A,T)\ge 2$ for any $k\ge 4$.
Thus $\Delta^k(A,T)=2$ or $3$ for $k\ge 4$.
The purpose of this paper is to determine this value.

\begin{theorem}\label{thm:main}
Let $M$ be a hyperbolic $3$-manifold with a torus boundary component $T_0$ and 
suppose that there are two slopes $\alpha$, $\beta$ on $T_0$ such that $M(\alpha)$ is annular
and $M(\beta)$ is toroidal.
If $\Delta(\alpha,\beta)=3$, then $\partial M$ is a union of at most three tori.
In particular, $\Delta(A,T)^k=2$ for any $k\ge 4$.
\end{theorem}

This gives a partial answer to \cite[Question 5.3]{Go4}.

In Section \ref{sec:pre}, we prepare the basic facts about labelled graphs.
In particular, the key is Lemma \ref{lem:key} which claims that neither graph contains both a black Scharlemann cycle
and a white Scharlemann cycle.
Section \ref{sec:one} is devoted to a special case where one graph has a single vertex, and
Section \ref{sec:two} deals with the case where the graph on the annulus has two vertices.
Section \ref{sec:generic} deals with the generic case.
To prove Theorem \ref{thm:main}, we need to consider the situation that $M(\beta)$ contains a Klein bottle.
This case is treated in Sections \ref{sec:klein} and \ref{sec:p=2}.

\section{Preliminaries}\label{sec:pre}

An annulus or torus is \textit{essential\/} if
it is incompressible, boundary-incompressible and is not boundary-parallel.
For two slopes $\alpha$ and $\beta$ on $T_0$, we suppose that $M(\alpha)$ is annular
and $M(\beta)$ is toroidal.
That is, $M(\alpha)$ (resp.\ $M(\beta)$) contains an essential annulus (resp.\ torus).

To prove Theorem \ref{thm:main}, we assume that $\Delta(\alpha,\beta)=3$
and $\partial M$ is not a union of at most three tori for contradiction, throughout the paper.

\begin{lemma}\label{lem:irr}
$M(\alpha)$ and $M(\beta)$ are irreducible and boundary-irreducible.
\end{lemma}

\begin{proof}
Since $M$ is large in the sense of \cite{W2}, $M(\alpha)$ and $M(\beta)$ are irreducible by \cite[Theorems 4.1 and 5.1]{W2}.
Boundary-irreducibility follows from \cite{GL3,GW2}.
\end{proof}

Let $\widehat{S}$ be an essential annulus in $M(\alpha)$.
For a core $K_\alpha$ of the attached solid torus $V_\alpha$, we can assume that
$\widehat{S}$ meets $K_\alpha$ transversely.
Then $\widehat{S}\cap V_\alpha$ is a disjoint union of meridian disks of $V_\alpha$, $u_1$, $u_2$, \dots, $u_s$,
numbered successively along $V_\alpha$, and $s$ can be chosen to be minimal among all essential annuli.
Similarly, we consider an essential torus $\widehat{T}$ in $M(\beta)$,
meeting a core $K_\beta$ of $V_\beta$ transversely.
Then $\widehat{T}\cap V_\beta$ is a union of meridian disks $v_1,v_2,\dots,v_t$, and $t$ is chosen to be minimal.
Let $S=\widehat{S}\cap M$ and $T=\widehat{T}\cap M$.
We can assume that no circle component in $S\cap T$ bounds a disk in $S$ or $T$, since both surfaces are incompressible.

In the usual way (\cite{CGLS,Go2,GW}), the arc components of $S\cap T$ define labelled graphs $G_S$ on $\widehat{S}$ and $G_T$ on $\widehat{T}$.
The vertices of $G_S$ (resp.\ $G_T$) are $u_1,u_2,\dots,u_s$ (resp.\ $v_1,v_2,\dots,v_t$).
For an edge of $G_S$, if its endpoint lies in $\partial u_i\cap \partial v_j$, then the point is labelled $j$ at $u_i$.
Thus the sequence of labels $1,2,\dots,t$ is repeated three times around each $u_i$, and so $u_i$ has degree $3t$.
Similarly, the edges of $G_T$ are labelled, and the sequence $1,2,\dots,s$ appears three times around $v_j$.
An edge with label $i$ at one of its endpoints is called an \textit{$i$-edge}.
Also, an edge with labels $i$ and $j$ is called a \textit{$\{i,j\}$-edge}.
An edge is said to be \textit{level\/} if its endpoints have the same label.
Notice that there is one-one correspondence between the edges of $G_S$ and $G_T$, and that
neither graph contains a \textit{trivial loop\/}, which bounds a $1$-sided disk face.
Throughout the paper, two graphs on a surface are considered to be equivalent if there is a homeomorphism of the surface sending one graph
to the other.

Each vertex of $G_S$ is given a sign, according to the sign of the intersection point of
$K_\alpha$ with $\widehat{S}$ with respect to some chosen orientations of $M$, $\widehat{S}$ and $K_\alpha$.
Similarly, we give a sign to each vertex of $G_T$.
Two vertices are \textit{parallel\/} if they have the same sign, otherwise they are \textit{antiparallel}.
An edge is \textit{positive\/} if it connects parallel vertices.
Otherwise, it is \textit{negative}.
In particular, a loop is positive.
A point at a vertex is called a \textit{positive edge endpoint\/} if there is a positive edge incident to there.
Otherwise, it is a \textit{negative edge endpoint}.

For a graph $G=G_S$ or $G_T$, let $G^+$ denote the subgraph consisting of all vertices and all positive edges of $G$.
Also, $G^+_x$ be the subgraph of $G^+$ consisting of all vertices and all $x$-edges of $G^+$ for a label $x$.
A disk face of $G^+_x$ is called an \textit{$x$-face}.
The \textit{reduced graph\/} $\overline{G}$ of $G$ is obtained from $G$ by amalgamating each family of mutually parallel edges into a single edge.

A cycle $\sigma$ consisting of positive edges is a \textit{Scharlemann cycle\/} if it bounds a disk face of the graph, and all the edges in $\sigma$ have
the same pair of labels $\{i,i+1\}$ at their endpoints, called the \textit{label pair\/} of $\sigma$.
The \textit{length\/} of $\sigma$ is the number of edges in $\sigma$.
In particular, a Scharlemann cycle of length two is called an \textit{$S$-cycle}.
If $\sigma$ is surrounded by a cycle $\tau$, that is, each edge of $\tau$ is immediately parallel to an edge of $\sigma$,
then $\tau$ is called an \textit{extended Scharlemann cycle} (see \cite{GL2}).

\begin{lemma}\label{lem:common}
\begin{itemize}
\item[(1)] There are no two edges which are parallel in both graphs.
\item[(2)] \textup{(The parity rule)} An edge is positive in one graph if and only if it is negative in the other.
\item[(3)] The edges of a Scharlemann cycle in $G_S$ \textup{(}resp.\ $G_T$\textup{)} do not lie in a disk in $\widehat{T}$ \textup{(}resp.\ $\widehat{S}$\textup{)}.
\item[(4)] If $G_S$ \textup{(}resp.\ $G_T$\textup{)} contains a Scharlemann cycle, then $\widehat{T}$ \textup{(}resp.\ $\widehat{S}$\textup{)}
is separating, and so $t$ \textup{(}resp.\ $s$\textup{)} is even.
\item[(5)] If $t>2$ \textup{(}resp. $s>2$\textup{)}, then $G_S$ \textup{(}resp.\ $G_T$\textup{)} cannot contain an extended Scharlemann cycle.
\end{itemize}
\end{lemma}

\begin{proof}
(1) This is \cite[Lemma 2.1]{Go2}.  See also \cite[Lemma 2.2]{GW}.
(2) can be found in \cite[p.279]{CGLS}.
(3) and (4) are \cite[Lemma 2.2]{GW}.
(5) For $G_S$, this is \cite[Theorem 3.2]{GL2}. 
For $G_T$, we refer \cite{LOT}.  Remark that only extended $S$-cycles are considered in \cite{GW,W2}.
\end{proof}


\begin{theorem}\label{thm:noklein}
$M(\beta)$ does not contain a Klein bottle meeting a core $K_\beta$ of $V_\beta$ in at most $t/2$ points.
\end{theorem}

This will be proved in Sections \ref{sec:klein} and \ref{sec:p=2}.

If $G_S$ contains a Scharlemann cycle, $M(\beta)$ is split into two pieces $\mathcal{B}$ and $\mathcal{W}$ along $\widehat{T}$.
We call them the \textit{black side\/} and the \textit{white side\/} of $\widehat{T}$, respectively.
Also, a disk face of $G_S$ is said to be \textit{black\/} or \textit{white\/}, according as it lies in $\mathcal{B}$ or $\mathcal{W}$.
In particular, a Scharlemann cycle whose disk face is black (white) is called a \textit{black \textup{(}white\textup{)} Scharlemann cycle}.
This is similar for $G_T$.

For the remainder of the section, let $H_{i,i+1}$ be that part of $V_\beta$ between $v_i$ and $v_{i+1}$.

\begin{lemma}\label{lem:key}
Neither graph contains a black Scharlemann cycle and a white Scharlemann cycle simultaneously.
\end{lemma}

\begin{proof}
Assume that $G_S$ contains a black Scharlemann cycle $\sigma_1$ and a white Scharlemann cycle $\sigma_2$.
Let $D_i$ be the disk face bounded by $\sigma_i$, and $\{k_i,k_{i}+1\}$ be the label pair of $\sigma_i$.
Let $X=N(\widehat{T}\cup H_{k_1,k_1+1}\cup H_{k_2,k_2+1}\cup D_1\cup D_2)$.
Then $\partial X$ consists of two tori $T_1$ and $T_2$, each of which intersects $K_\beta$ fewer than $t$ times.
By the minimality of $\widehat{T}$, each $T_i$ is boundary-parallel or bounds a solid torus in $M(\beta)$.
Thus $\partial M$ consists of at most three tori, contradicting our assumption.

A similar construction works for $G_T$.
By using the disk faces bounded by a black Scharlemann cycle and a white Scharlemann cycle in $G_T$,
we obtain two annuli $S_1$ and $S_2$, each of which intersects $K_\alpha$ fewer than $s$ times.
By the minimality of $\widehat{S}$, each $S_i$ is boundary parallel.
Hence $\partial M$ is a union of two tori, a contradiction.
\end{proof}


\begin{lemma}\label{lem:parallel-max}
$G_S$ satisfies the following.
\begin{itemize}
\item[(1)] There are no two $S$-cycles with disjoint label pairs.
\item[(2)] Any family of mutually parallel positive edges in $G_S$ contains at most $t/2+1$ edges.
If $\widehat{T}$ is non-separating, then it contains at most $t/2$ edges.
\item[(3)] Either any family of mutually parallel negative edges in $G_S$ contains at most $t$ edges,
or all vertices of $G_T$ are parallel.
\end{itemize}
\end{lemma}

\begin{proof}
(1) Let $\sigma_1$ and $\sigma_2$ be $S$-cycles with disjoint label pairs.
Let $\{k_i,k_i+1\}$ be the label pair of $\sigma_i$, and let $D_i$ be the face bounded by $\sigma_i$, $i=1,2$.
Shrinking $H_{k_i,k_i+1}$ to its core in $H_{k_i,k_i+1}\cup D_i$ gives a M\"{o}bius band $B_i$ whose boundary is
essential on $\widehat{T}$ by Lemma \ref{lem:common}(3).
By $\partial B_1$ and $\partial B_2$, $\widehat{T}$ is split into two annuli $A_1$ and $A_2$.
If $A_i$ contains $a_i$ vertices in its interior, then the Klein bottle $F_i=B_1\cup A_i\cup B_2$ meets
$K_\beta$ in $a_i+2$ points.
The torus $\partial N(F_i)$ is incompressible by the irreducibility of $M(\beta)$.
If it is boundary parallel in $M(\beta)$, then $M(\beta)=N(F_i)$ is not toroidal.
Hence this torus is essential, and so $2(a_i+2)\ge t$, giving $a_i\ge t/2-2$.
Since $a_1+a_2=t-4$, $a_1=a_2=t/2-2$.
But this contradicts Theorem \ref{thm:noklein}. 

(2) If $t>2$, then such family contains at most $t/2+2$ edges by \cite[Lemma 1.4]{W}, and moreover,
if it contains $t/2+2$ edges, then it contains two $S$-cycles with disjoint label pairs, which contradicts (1).
If $t=2$, then three parallel edges contain a black $S$-cycle and a white $S$-cycle, contradicting Lemma \ref{lem:key}.
When $t=1$, $G_T$ contains only positive edges, and so $G_S$ has no positive edges by the parity rule.

If a family of parallel positive edges contains more than $t/2$ edges, then the family contains an $S$-cycle.
The second claim follows immediately from Lemma \ref{lem:common}(4).

(3) See \cite[Lemma 2.3]{GW} for $t>2$. 
Assume $t=2$ and that two vertices of $G_T$ are antiparallel.
Suppose that $G_S$ has three mutually parallel negative edges $e_1$, $e_2$, $e_3$, numbered successively.
Then all are level by the parity rule.
Hence we can assume that $e_1$ and $e_3$ have label $1$.
Since two loops at each vertex of $G_T$ are parallel (see \cite[Lemma 5.2]{Go2}),
$e_1$ and $e_3$ correspond to parallel loops at $v_1$ in $G_T$.
This contradicts Lemma \ref{lem:common}(1).
\end{proof}

\begin{lemma}\label{lem:GT}
$G_T$ satisfies the following.
\begin{itemize}
\item[(1)] If $s>2$, then any family of mutually parallel positive edges in $G_T$ contains at most $s/2+1$ edges.
If $\widehat{S}$ is non-separating, then it contains at most $s/2$ edges.
\item[(2)] Any family of mutually parallel negative edges in $G_T$ contains at most $s$ edges.
\item[(3)] All Scharlemann cycles in $G_T$ have the same label pair.
\end{itemize}
\end{lemma}

\begin{proof}
These are \cite[Lemma 2.5]{GW}.
\end{proof}

\begin{lemma}\label{lem:GS}
$G_S$ satisfies the following.
\begin{itemize}
\item[(1)] At most two labels can be labels of $S$-cycles.
\item[(2)] At most four labels can be labels of Scharlemann cycles.
\end{itemize}
\end{lemma}

\begin{proof}
(1) If there are three labels of $S$-cycles, there are two $S$-cycles with disjoint label pairs and
with the same color by Lemma \ref{lem:key}.
This is impossible by Lemma \ref{lem:parallel-max}(1).

(2)
If not, $G_S$ has three Scharlemann cycles $\sigma_1$, $\sigma_2$, $\sigma_3$ with mutually
disjoint label pairs and with the same color Lemma \ref{lem:key}.
Let $D_i$ be the face bounded by $\sigma_i$, and let $\{k_i,k_i+1\}$ be the label pair of $\sigma_i$.
We can assume that $D_i\subset \mathcal{B}$.
On $\widehat{T}$, there are mutually disjoint annuli $A_i$ which contain the edges of $\sigma_i$ respectively.
Define $M_i=N(A_i\cup H_{k_i,k_i+1}\cup D_i)\subset \mathcal{B}$.
Let $B_i=\mathrm{cl}\,(\partial M_i-A_i)$.
Then a new torus $T_i=(\widehat{T}-A_i)\cup B_i$ meets $K_\beta$ fewer than $t$ times.
Hence $T_i$ is compressible or boundary parallel.
If one of $T_i$ is compressible, the argument in the proof of \cite[Theorem 3.5]{GL2} without any change gives a contradiction.
Thus any $T_i$ is boundary parallel.
Let $Z_1=\mathrm{cl}\,(\mathcal{B}-M_1)$.
Then $Z_1=T^2\times I$.
Since $M_2\subset Z_1$, $M_2$ is a solid torus, and moreover, $B_2$ is parallel to $A_2$ through $M_2$.
This contradicts \cite[Claim 3.6]{GL2}.
\end{proof}

If $G_S$ contains a Scharlemann cycle with label $i$, then $i$ is called an \textit{$S$-label}.
Otherwise, $i$ is called a non-$S$-label.

\begin{lemma}\label{lem:key2}
Let $t\ge 3$.
Any $x$-face in $G_S$ has at least $4$ sides for a non-$S$-label $x$.
\end{lemma}

\begin{proof}
Assume not.
By Lemmas \ref{lem:common}(5) and \ref{lem:parallel-max}(2), $G_S$ cannot contain a two-sided $x$-face. 
Let $D$ be a $3$-sided $x$-face in $G_S$.
By \cite[Proposition 5.1]{HM}, $D$ contains a Scharlemann cycle.
Since $G_S$ cannot contain an extended Scharlemann cycle by Lemma \ref{lem:common}(5), $D$ contains an $S$-cycle.
By using Lemma \ref{lem:GS}(1), the proof of \cite[Lemma 5.1]{GL2} shows that $D$ contains an $S$-cycle with face $f$,
and the bigon $g_1$ and the $3$-gon $g_2$ adjacent to $f$ have only two kinds of corners.
See \cite[Figure 5.4]{GL2}.
For convenience, we assume that $f$ has two $(1,2)$-corners, and $g_i$ has $(t,1)$- and $(2,3)$-corners.
Let $A_{i,i+1}$ be the annulus in $\partial V_\beta$ between $v_i$ and $v_{i+1}$.
Then $(\widehat{T}-\mathrm{Int}\,(v_t\cup v_1\cup v_2\cup v_3))\cup (A_{t,1}\cup A_{2,3})$ is a genus three closed surface,
on which $\partial g_1$ and $\partial g_2$ are homologically independent.
(This means that the genus three closed surface will be compressed to a torus along $g_1$ and $g_2$.)
Hence $N(\widehat{T}\cup H\cup f\cup g_1\cup g_2)$ has two torus boundary components, where
$H$ is the part of $V_\beta$ between $v_t$ and $v_3$, containing $v_1$.
Since each torus meets $K_\beta$ fewer than $t$, they are inessential in $M(\beta)$.
Then $M$ is bounded by at most three tori as in the proof of Lemma \ref{lem:key}, a contradiction.
\end{proof}

\section{Case where one graph has a single vertex}\label{sec:one}

In this section, we treat the case where $s=1$ or $t=1$.

\begin{lemma}
$s\ne 1$.
\end{lemma}

\begin{proof}
Assume $s=1$.
Since the vertex $u_1$ of $G_S$ has degree $3t$, $t$ must be even.
Also, all edges of $G_S$ are positive and parallel.
If $t>2$, then $3t/2\le t/2+1$, giving $t\le 1$, a contradiction.
Hence $t=2$.
But then $G_S$ contains a black $S$-cycle and a white $S$-cycle, which contradicts Lemma \ref{lem:key}. 
\end{proof}

\begin{lemma}\label{lem:t1}
If $t=1$, then $s=2$.
\end{lemma}

\begin{proof}
Assume $t=1$ and $s\ge 3$.
There are $3s/2$ edges (so, $s$ is even) in $G_T$, which are divided into at most three families of mutually parallel edges (see \cite[Lemma 5.1]{Go2}).
Since each family contains at most $s/2+1$ edges by Lemma \ref{lem:GT}(1), $G_T$ has at least two families.
If there are only two families, then $4(s/2+1)\ge 3s$ gives $s\le 4$.
Hence $s=4$.
Then $G_T$ consists of two families of three mutually parallel edges.
By examining the labels, this contradicts the parity rule.
Therefore, $G_T$ contains three families.

We denote by $G_T\cong H(q_1,q_2,q_3)$ when each family contains
$q_1,q_2,q_3$ edges, respectively.
Note that $H(q_1,q_2,q_3)$ is invariant under any permutation of the $q_i$'s.
If $q_i\le s/2$ for any $i$, then $6\cdot s/2\ge 2(q_1+q_2+q_3)=3s$ gives
$q_1=q_2=q_3=s/2$, and so $G_T\cong H(s/2,s/2,s/2)$.
It is easy to see that $G_T$ contains an extended Scharlemann cycle of length three, which is impossible by Lemma \ref{lem:common}(5).
Hence we may assume $q_1=s/2+1$.
Let $Q_i$ denote the family of parallel edges containing $q_i$ edges.
Then $Q_1$ contains an $S$-cycle at one end (see \cite[Lemma 1.4]{W}).
By examining the labels, $q_2+q_3\equiv 0 \pmod{s}$ and $q_2+q_3\equiv s-2 \pmod{s}$.
This implies $s=2$, a contradiction.
\end{proof}

\begin{lemma}
$t\ne 1$.
\end{lemma}

\begin{proof}
Assume $t=1$. By Lemma \ref{lem:t1}, $s=2$.
Then $G_T\cong H(3,0,0)$, $H(2,1,0)$ or $H(1,1,1)$.
If $G_T\cong H(3,0,0)$ or $H(1,1,1)$, then $G_T$ contains a black Scharlemann cycle and a white Scharlemann cycle, contradicting Lemma \ref{lem:key}.
Clearly, $H(2,1,0)$ contradicts the parity rule.
\end{proof}

\section{Case where $s=2$}\label{sec:two}

\begin{figure}[tb]
\includegraphics*[scale=0.45]{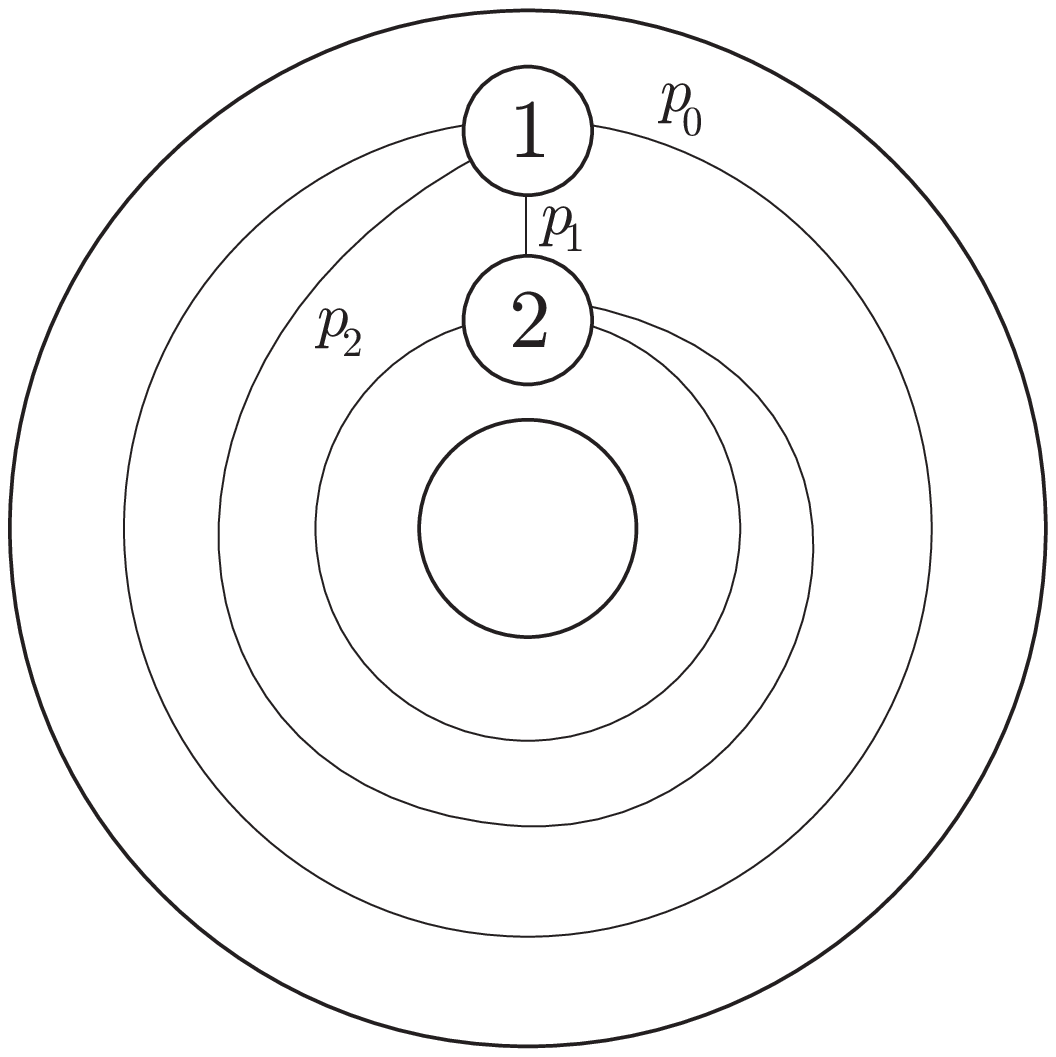}
\caption{}\label{fig:annulus}
\end{figure}

In this section, we consider the case where $s=2$ and $t\ge 2$.
Then the reduced graph $\overline{G}_S$ of $G_S$ is a subgraph of the graph shown in Figure \ref{fig:annulus}.
Notice that $u_1$ and $u_2$ are incident to the same number of loops in $G_S$.
We denote by $G_S\cong G(p_0,p_1,p_2)$ when $u_i$ is incident to $p_0$ loops, and
the other two families of parallel edges contain $p_1$ and $p_2$ edges, respectively.
Clearly, $G(p_0,p_1,p_2)$ is equivalent to $G(p_0,p_2,p_1)$.
We divide the argument into two cases.

\subsection{The two vertices of $G_S$ are parallel}\label{subsec:parallel}

\begin{lemma}
$t>2$.
\end{lemma}

\begin{proof}
Assume $t=2$.
Then $p_i\le 2$ for any $i$ by Lemma \ref{lem:parallel-max}(2), and so $G_S\cong G(2,2,0)$, $G(2,1,1)$ or $G(1,2,2)$.
In any case, $G_S$ contains a black Scharlemann cycle and a white Scharlemann cycle, since any disk face is a Scharlemann cycle.
This contradicts Lemma \ref{lem:key}.
\end{proof}

\begin{lemma}
$\widehat{T}$ is separating and $t=4$.
\end{lemma}

\begin{proof}
Since $3t=2p_0+p_1+p_2$, $p_i>t/2$ for some $i$.
Thus $G_S$ contains an $S$-cycle, and so $\widehat{T}$ is separating and $t$ is even by Lemma \ref{lem:common}(4).
Hence we have $t\ge 4$.
Notice that $p_i\le t/2+1$ for any $i$ by Lemma \ref{lem:parallel-max}(2).
Hence $3t\le 4(t/2+1)=2t+4$, giving $t\le 4$.
\end{proof}

\begin{proposition}
The two vertices of $G_S$ cannot be parallel.
\end{proposition}

\begin{proof}
Since $p_i\le 3$ for any $i$ by Lemma \ref{lem:parallel-max}(2),
$G_S\cong G(3,3,3)$.
Once we fix labeling around $u_1$, there are only two possibilities for labelling around $u_2$
by the parity rule.
In any case, $G_S$ contains a black Scharlemann cycle and a white Scharlemann cycle, a contradiction.
\end{proof}

\subsection{The two vertices of $G_S$ are antiparallel}\label{subsec:antiparallel}

\begin{lemma}
$p_0\ne 0$.
\end{lemma}

\begin{proof}
If $p_0=0$, then all edges of $G_S$ connect $u_1$ with $u_2$.
Hence all edges of $G_T$ are positive by the parity rule.
Notice that any disk face of $G_T$ is a Scharlemann cycle.
Since $G_T$ has $3t$ edges, it contains at least $2t$ disk faces.
If these disk faces have the same color, then $G_T$ has at least $4t$ edges, a contradiction.
Thus $G_T$ contains a black Scharlemann cycle and a white Scharlemann cycle, contradicting Lemma \ref{lem:key}.
\end{proof}

\begin{lemma}\label{lem:s2t3}
$G_S\cong G(t/2,t,t)$, $G(t/2+1,t,t-2)$ or $G(t/2+1,t-1,t-1)$.
\end{lemma}

\begin{proof}
Since $p_0\ne 0$, $G_S$ contains a positive edge, and not all the vertices of $G_T$ are parallel.
This implies $p_i\le t$ for $i=1,2$ by Lemma \ref{lem:parallel-max}(3).
By $3t=2p_0+p_1+p_2\le 2p_0+2t$, $p_0\ge t/2$.
If $p_0=t/2$, then $G_S\cong G(t/2,t,t)$.
If $p_0>t/2$, then $G_S$ contains an $S$-cycle, and hence $\widehat{T}$ is separating and $t$ is even.
Thus $p_0=t/2+1$, and then the conclusion follows immediately.
\end{proof}

\begin{lemma}\label{lem:noS}
$G_T$ cannot contain an $S$-cycle.
\end{lemma}

\begin{proof}
Let $\sigma$ be an $S$-cycle in $G_T$ whose disk face is $f$.
The edges of $\sigma$ form an essential cycle in $\widehat{S}$ by Lemma \ref{lem:common}(3).
Let $H$ be the part of $V_\alpha$ between $u_1$ and $u_2$ meeting $\partial f$.
Then shrinking $H$ into its core in $H\cup f$ gives a M\"{o}bius band $B'$ whose boundary is an essential loop on $\widehat{S}$.
The union of $B'$ and an annulus between $\partial B'$ and $\partial S$ gives a M\"{o}bius band $\widehat{B}$
properly embedded in $M(\beta)$ which meets $K_\alpha$ in one point.
Let $X=N(\widehat{B})$ and let $W=M(\alpha)-\mathrm{Int}\,X$.
Then the frontier $\widehat{Q}$ of $X$ is an incompressible annulus.
If $\widehat{Q}$ is boundary parallel, then $M(\alpha)$ has a single torus boundary, a contradiction.
Hence $\widehat{Q}$ is essential.
Let $Q=\widehat{Q}\cap M$, and let $A=\partial V_\alpha\cap W$.
Then $F=Q\cup A$ is a twice-punctured torus.

Let $B=\widehat{B}\cap M$.
If $B$ is compressible in $M$, then let $\delta$ be a compressing disk for $B$.
Since $\partial \delta$ is orientation-preserving on $B$, it bounds a disk in $\widehat{B}$ or is parallel to $\partial \widehat{B}$.
The former implies that $M$ contains a properly embedded M\"{o}bius band, contradicting the hyperbolicity of $M$.
The latter means that $M(\alpha)$ contains a projective plane, and so $M(\alpha)$ is reducible, contradicting Lemma \ref{lem:irr}.
Hence $B$ is incompressible.
Also, if $B$ is boundary compressible, then $K_\alpha$ can be isotoped to the core of $\widehat{B}$ by using a boundary compressing disk.
Then $M$ contains an essential annulus.
Hence $B$ is boundary incompressible.

We construct another graph pair $\{G_B,G_T^B\}$ from $B$ and $T$ in the usual way.
There is no trivial loop in each graph.
Note that $G_B$ has a single vertex, and $G_T^B$ consists of $t$ vertices of degree three and $3t/2$ edges.
In fact, the double cover of $G_B$ is a subgraph of the graph shown in Figure \ref{fig:annulus}.
By an Euler characteristic calculation, $G_T^B$ contains a disk face $D'$.
Let $D=D'\cap W$.
Notice that $\partial D$ is essential on $F$.
For, $\partial D$ runs on $Q$ and $A$ alternatively,
and $\partial D\cap Q$ consists of arcs as shown in Figure \ref{fig:annulus}.
Surgering $F$ along $D$ gives either an annulus or a disjoint union of an annulus and a torus, according as $\partial D$ is non-separating or
separating on $F$.
In any case, the resulting surface is disjoint from $K_\alpha$.
Hence the annulus component is boundary parallel, and the torus component, if it exists, is inessential.
Thus $M(\alpha)$ is bounded by at most two tori, a contradiction.
\end{proof}

\begin{lemma}\label{lem:t2}
$t=2$.
\end{lemma}

\begin{proof}
Assume $t>2$.  There are three possibilities for $G_S$ by Lemma \ref{lem:s2t3}.

If $G_S\cong G(t/2,t,t)$, then $G_T$ has $2t$ positive edges by the parity rule.
Hence $G_T^+$ has at least $t$ disk faces.
Notice that such disk face is also a face of $G_T$, and so it is bounded by a Scharlemann cycle.
Hence we may assume that such disk faces are all black by Lemma \ref{lem:key}.
Also, such disk face has at least three sides by Lemma \ref{lem:noS}.
Thus there are at least $3t$ positive edges, a contradiction.

If $G_S\cong G(t/2+1,t,t-2)$, then the same argument yields a contradiction, unless $t=4$.
(Notice that $p_0=t/2+1$ implies that $\widehat{T}$ is separating and $t$ is even by Lemma \ref{lem:common}(4).)
Suppose $t=4$ and $G_S\cong G(3,4,2)$.
Let $Q$ be the family of $4$ negative edges in $G_S$, and let $\sigma$ be the associated permutation to $Q$.
That is, each edge of $Q$ has label $i$ at $u_1$ and $\sigma(i)$ at $u_2$.
If $\sigma$ is the identity, then $G_S$ contains two $S$-cycles with disjoint label pairs, which contradicts Lemma \ref{lem:parallel-max}.
Hence $\sigma=(13)(24)$.
In this case, $G_T$ is uniquely determined.
First, the edges of two $S$-cycles with label pair $\{3,4\}$ form essential cycles.
The edges of $Q$ form two essential cycles by Lemma \ref{lem:common}(1).
By examining labels, two edges between $v_1$ and $v_2$ turn out to be parallel.
See Figure \ref{fig:4-7}.

\begin{figure}[tb]
\includegraphics*[scale=0.5]{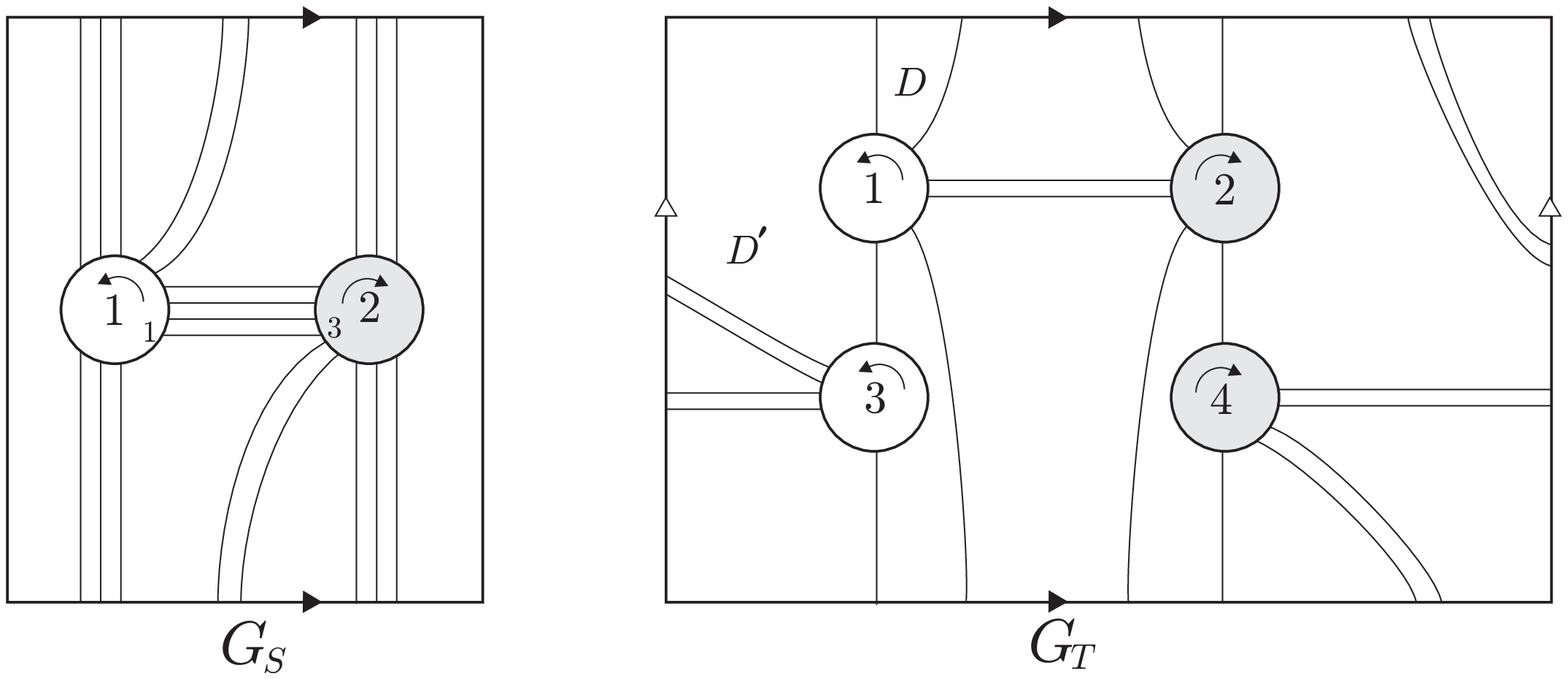}
\caption{}\label{fig:4-7}
\end{figure}

Then $G_T$ has a Scharlemann cycle of length three with face $D$.
Thus $M(\alpha)$ is split into two pieces $\mathcal{B}$ and $\mathcal{W}$ along $\widehat{S}$.
We may assume that $D\subset \mathcal{B}$.
Let $H=V_\alpha\cap \mathcal{B}$.
Let $X=N(\widehat{S}\cup H\cup D)\subset \mathcal{B}$.
By the minimality of $\widehat{S}$, the annulus $\mathrm{cl}(\partial X-\widehat{S})$ is boundary parallel.
Thus $\partial \mathcal{B}$ is a torus.
Let $D'$ be the white face as shown in Figure \ref{fig:4-7}.
Similarly, we can see that $\partial \mathcal{W}$ is a torus by using $D'$.
Thus $M(\alpha)$ is bounded by a single torus, a contradiction.

If $G_S\cong G(t/2+1,t-1,t-1)$, then two families of loops at $u_1$ and $u_2$ contain
$S$-cycles.  Hence $t$ is even.
By examining labels, such $S$-cycle is located at one end of each family.
Then it is obvious that these two $S$-cycles have distinct colors, contradicting Lemma \ref{lem:key}.
\end{proof}

\begin{figure}[tb]
\includegraphics*[scale=0.3]{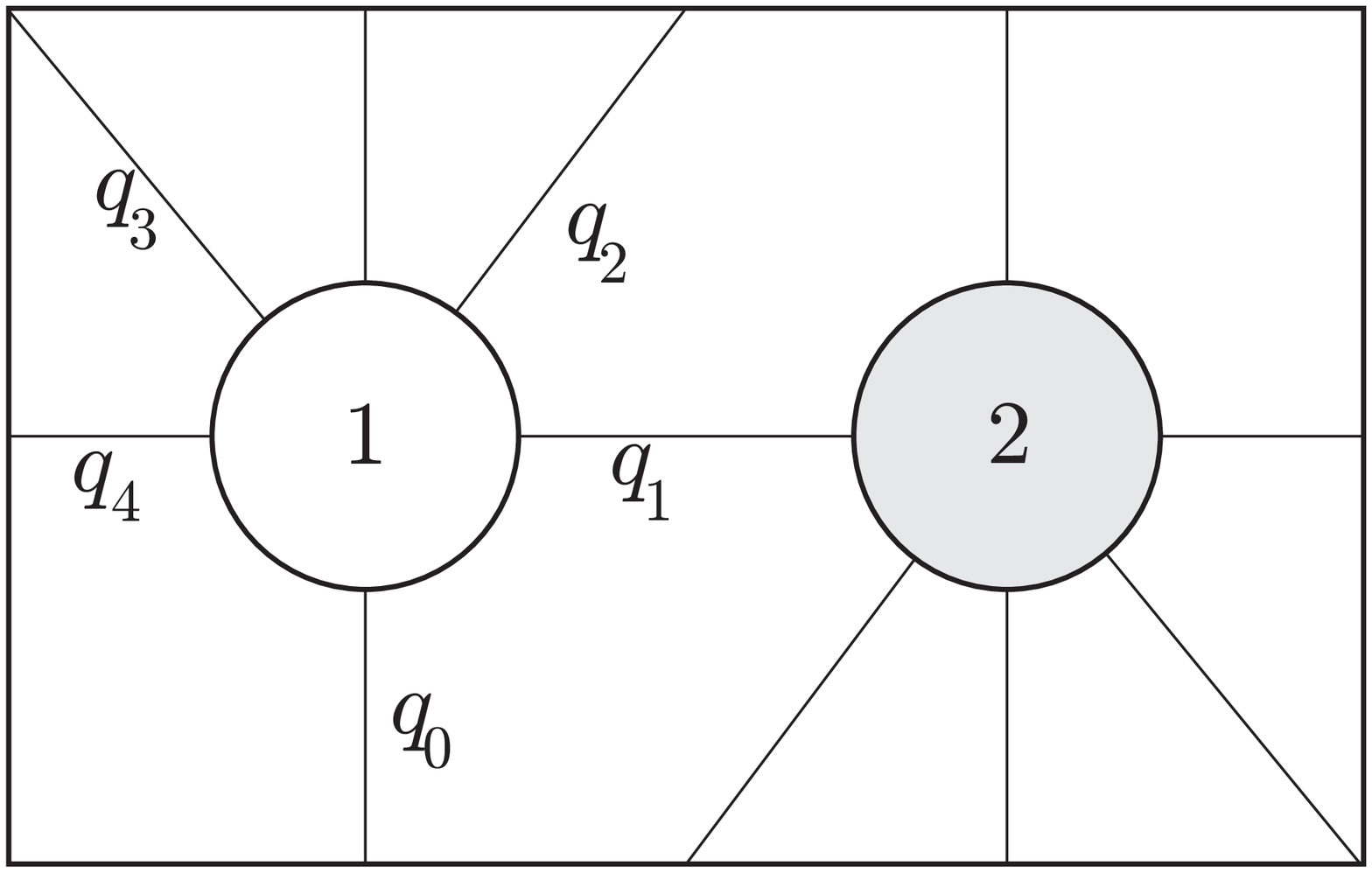}
\caption{}\label{fig:t2}
\end{figure}

By Lemma \ref{lem:t2}, $G_T$ has only two vertices.
The reduced graph $\overline{G}_T$ is a subgraph of the graph shown in Figure \ref{fig:t2} (see \cite[Lemma 5.2]{Go2}).
We say $G_T\cong H'(q_0,q_1,q_2,q_3,q_4)$, where $q_i$ denotes the number of edges in the family of parallel edges.
Note that
\begin{eqnarray*}
H'(q_0,q_1,q_2,q_3,q_4) &\cong & H'(q_0,q_3,q_4,q_1,q_2) \cong H'(q_0,q_4,q_3,q_2,q_1) \\
                        &\cong & H'(q_0,q_2,q_1,q_4,q_3).
\end{eqnarray*}

\begin{proposition}
The two vertices of $G_S$ cannot be antiparallel.
\end{proposition}

\begin{proof}
By Lemma \ref{lem:s2t3}, $G_S\cong G(1,2,2)$, $G(2,1,1)$ or $G(2,2,0)$.

If $G_S\cong G(1,2,2)$, then $G_T\cong H'(2,1,1,0,0)$ or $H'(2,2,0,0,0)$.
Then $G_T$ contains an $S$-cycle, contradicting Lemma \ref{lem:noS}.
If $G_S\cong G(2,1,1)$, then $G_S$ contains a black Scharlemann cycle and a white Scharlemann cycle, contradicting Lemma \ref{lem:key}.

Suppose $G_S\cong G(2,2,0)$.
Then $G_S$ contains two $S$-cycles $\rho_1$ and $\rho_2$ of the same color.
Let $f_i$ be its face for $i=1,2$, and let $A$ be the annulus part of $\partial V_\beta$ between $v_1$ and $v_2$, meeting $f_i$.
Notice that $q_0=1$ and $(q_1+q_2,q_3+q_4)=(3,1)$, $(2,2)$ or $(4,0)$, up to equivalence.
$(3,1)$ is impossible by the parity rule.
Thus $G_T\cong H'(1,1,1,2,0)$, $H'(1,2,0,2,0)$, $H'(1,1,1,1,1)$, or $H'(1,2,2,0,0)$.

First, $H'(1,1,1,2,0)$ contradicts the parity rule.
If $G_T\cong H'(1,2,0,2,0)$, $\partial f_1$ and $\partial f_2$ cannot be located 
on $T\cup A$ simultaneously.
Assume $G_T\cong H'(1,1,1,1,1)$.
Then there are two disjoint rectangles $R_1$ and $R_2$ in $A$ split by $\partial f_1\cup \partial f_2$ such that
$f_i\cup R_i$ gives a M\"{o}bius band $B_i$.
Thus we have two M\"{o}bius bands $B_1$ and $B_2$ whose boundaries are disjoint on $\widehat{T}$.
Hence $M(\beta)$ contains a Klein bottle as a union of $B_1$, $B_2$ and an annulus on $\widehat{T}$, meeting $K_\beta$ once.
This contradicts Theorem \ref{thm:noklein}.
Finally, assume $G_T\cong H'(1,2,2,0,0)$.
Then $G_T$ contains a $3$-gon $f$ and a bigon $g$.
Let $A'$ (resp.\ $A''$) be the part of $\partial V_\alpha$ between $u_1$ and $u_2$ meeting $\partial f$ (resp.\ $\partial g$).
Then $\partial f$ is a non-separating curve on the surface $S\cup A'$, so surgering $S\cup A'$ along $f$ gives rise to a boundary parallel annulus in $M(\alpha)$.
Thus $\widehat{S}$ is separating in $M(\alpha)$.
On the other hand, surgering $S\cup A''$ along $g$ gives rise to a surface disjoint from $K_\alpha$,
which is an annulus or a disjoint union of an annulus and a torus, according as $\partial g$ is non-separating or separating on $S\cup A''$.
As in the proof of Lemma \ref{lem:noS}, $M(\alpha)$ is bounded by at most two tori, a contradiction.
\end{proof}

\section{Generic case}\label{sec:generic}

Finally, we consider the case where $s\ge 3$ and $t\ge 2$.
Since all Scharlemann cycles of $G_T$ have the same label pair by Lemma \ref{lem:GT}(3), 
we can assume that $\{1,2\}$ is the label pair, if they exist.
Then these labels are $S$-labels of $G_T$, and
the vertices $u_1$ and $u_2$ are referred to as the \textit{$S$-vertices\/} of $G_S$.

\begin{lemma}\label{lem:xface}
$G_T$ does not contain an $x$-face for a non-$S$-label $x$.
\end{lemma}

\begin{proof}
This is Theorem 4.5 of \cite{LOT}.
\end{proof}

\begin{lemma}\label{lem:2t}
Any vertex of $G_S$, except $S$-vertices, has at least $2t$ positive edge endpoints.
\end{lemma}

\begin{proof}
Assume that $u_i$ is not an $S$-vertex. 
If it has at least $t+1$ negative edge endpoints,
then $G_T$ has at least $t+1$ positive $i$-edges.
Let $\Gamma_i$ be the subgraph of $G_T$ consisting of all vertices and all positive $i$-edges of $G_T$.
Then an Euler characteristic calculation shows that $\Gamma_i$ has a disk face, which is an $i$-face.
This contradicts Lemma \ref{lem:xface}.
\end{proof}

\begin{lemma}\label{lem:t}
An $S$-vertex of $G_S$, if it exists, has at least $t$ positive edge endpoints.
\end{lemma}

\begin{proof}
Let $u_1$ be an $S$-vertex.
Suppose that $u_1$ has $k$ negative edge endpoints.
Then $G_T$ has $k$ positive $1$-edges.
Hence $G_T$ has at least $k-t$ $1$-faces.
Recall that each $1$-face contains a Scharlemann cycle by \cite{HM}.
Thus there are at least $k-t$ Scharlemann cycles with label pair $\{1,2\}$.
Then there are at least $2(k-t)$ positive $1$-edges, since all Scharlemann cycles have the same color.
We have $2(k-t)\le k$, and so $k\le 2t$.
Hence $u_1$ has at least $t$ positive edge endpoints.
\end{proof}

Let us consider $G_S^+$, which consists of all vertices and all positive edges of $G_S$.
Let $\Lambda$ be a component of $G_S^+$.
If there is a disk $D$ in $\widehat{S}$ such that $\Lambda\subset \mathrm{Int}\,D$, then $\Lambda$ is said to have a \textit{disk support}.
Otherwise, there is an annulus $A$ in $\widehat{S}$, which is called an \textit{annulus support},
such that $\Lambda\subset \mathrm{Int}\,A$.
Clearly, the core of $A$ is parallel to the core of $\widehat{S}$.
Furthermore, if $\Lambda$ has a support $F$, which is a disk or an annulus, such that $F\cap G_S^+=\Lambda$, then
$\Lambda$ is called an \textit{extremal component\/} of $G_S^+$.
Clearly, if there is no component of $G_S^+$ with a disk support, then any component of $G_S^+$ is an extremal one with an annulus support.

Suppose that $\Lambda$ is an extremal component with support $F$.
A vertex $u$ is a \textit{cut vertex\/} if $\Lambda-u$ has more components than $\Lambda$.
We remark that $\Lambda$ may have loops.
Also, $u$ is called an \textit{interior vertex\/} if there is no arc $\xi$ in $F$ connecting $u$ to $\partial F$ such that $\xi\cap \Lambda=u$.
Otherwise, $u$ is called a \textit{boundary vertex}.
Furthermore, an \textit{interior edge\/} is an edge which cannot admit an arc $\xi$ connecting a middle point $x$ of the edge to $\partial F$
such that $\xi\cap \Lambda=x$.
The others are \textit{boundary edges\/}.
When $F$ is an annulus, a vertex $u$ is called a \textit{pinched vertex\/} if there is a spanning arc $\xi$ of $F$ such that
$\xi\cap \Lambda=u$, and a \textit{pinched edge\/} is defined similarly.
In particular, both endpoints of a pinched edge are pinched vertices.
Finally, $u$ is said to be \textit{good\/} if all positive edge endpoints at $u$ are successive.
Thus, if $u$ is neither a cut vertex nor a pinched vertex, then it is good.


A subgraph $B$ of $\Lambda$ is called a \textit{disk block\/} of $G_S^+$
if $B$ contains at most one cut vertex of $\Lambda$ and
there is a disk $D$ in $\widehat{S}$ such that $D\cap G_S^+=B$ and $\partial D\cap B$ is either empty or a single vertex.
We remark that a disk block is connected and that
a disk block cannot contain a loop which is essential in $\widehat{S}$, but it may contain
a loop which is inessential in $\widehat{S}$.
If $B$ has an $S$-vertex $u$, then $u$ must appear as a boundary vertex of $B$, because
the edges of a Scharlemann cycle in $G_T$ do not lie in a disk in $\widehat{S}$ by Lemma \ref{lem:common}(3).


\subsection{Case $t=2$}\label{subsec:t2}

To eliminate the case where $t=2$, we prove three lemmas.
Recall that any non-$S$-vertex has at least $4$ positive edge endpoints by Lemma \ref{lem:2t},
while any $S$-vertex has at least two positive edge endpoints by Lemma \ref{lem:t}.

\begin{lemma}\label{lem:t2pre}
Any component of $G_S^+$ has an annulus support, and hence is extremal.
\end{lemma}

\begin{proof}
If $G_S^+$ has a component with a disk support, then there is an extremal component $\Lambda$ with a disk support.
By Lemmas \ref{lem:2t} and \ref{lem:t}, it contains at least two vertices.
Any vertex of $\Lambda$, except a cut vertex, is good, and has at least $4$ successive positive edge endpoints by Lemma \ref{lem:2t}.
Hence $\Lambda$ has a black face and a white face, which contradicts Lemma \ref{lem:key}, because
any disk face of $G_S^+$ is bounded by a Scharlemann cycle.

Therefore we have shown that any component of $G_S^+$ has an annulus support.
Also, this implies that any component is extremal.
\end{proof}

\begin{lemma}\label{lem:t2dblock}
$G_S^+$ has at most two disk blocks, each of which consists of two vertices and a pair of parallel edges.
In particular, a non-cut vertex is an $S$-vertex.
\end{lemma}

\begin{proof}
Let $B$ be a disk block.
If $B$ has an interior edge, then there is a black face and a white face, contradicting Lemma \ref{lem:key}.
Hence $B$ has no interior edge.
Thus $B$ is either a single edge or a cycle.
However, the former is impossible by Lemmas \ref{lem:2t} and \ref{lem:t}.
Hence $B$ is a cycle.
If the length of $B$ is more than two,
then there is a non-cut vertex, which is not an $S$-vertex, contradicting Lemma \ref{lem:2t}.
Hence $B$ is length two, and Lemma \ref{lem:2t}
implies that a non-cut vertex must be an $S$-vertex.

Since $G_S$ has at most two $S$-vertices, there are at most two disk blocks.
\end{proof}

\begin{lemma}\label{lem:t2pre2}
Any component of $G_S^+$ containing a non-$S$-vertex is a cycle of bigons.
\end{lemma}

\begin{proof}
Let $\Lambda$ be a component containing a non-$S$-vertex $u$.
Recall that every face of $\Lambda$ is a disk bounded by a Scharlemann cycle.
Hence $\Lambda$ has no interior vertex.

First, assume that $\Lambda$ has no cut vertex.
Recall that any non-$S$-vertex has at least $4$ positive edge endpoints.
Also, $\Lambda$ has at most one $S$-vertex.
If a non-$S$-vertex is not pinched, then $\Lambda$ has a black face and a white face.
Hence any non-$S$-vertex is pinched, and has degree $4$.
Thus $\Lambda$ is either a cycle of bigons, or a cycle of bigon added one bivalent vertex, which is an $S$-vertex.
See Figure \ref{fig:bigon-cycle}.

\begin{figure}[tb]
\includegraphics*[scale=0.6]{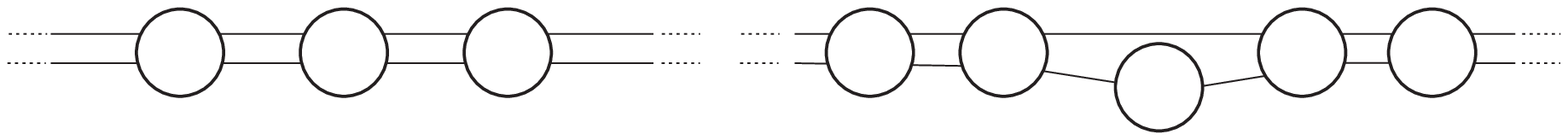}
\caption{}\label{fig:bigon-cycle}
\end{figure}

Suppose that $\Lambda$ contains a bivalent $S$-vertex $u_1$, say.
Then the configuration of $G_S$ near $u_1$ looks like Figure \ref{fig:bigon-S}.
Notice that $u_1$ has $4$ negative edges, so $G_T$ has at least two $1$-faces, which must be bigons bounded by $S$-cycles.
Thus $G_T$ contains two $S$-cycles.
Let $D$ be the disk face as shown there.

\begin{figure}[tb]
\includegraphics*[scale=0.6]{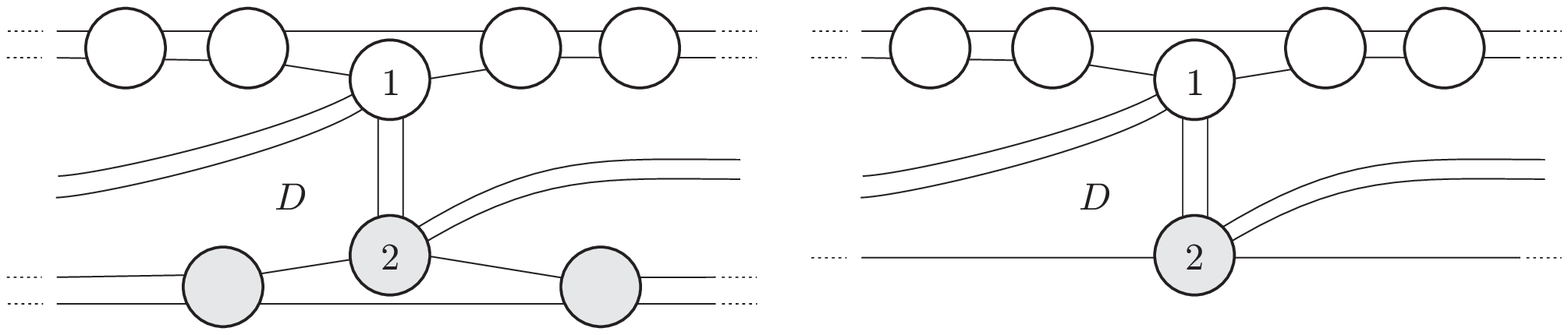}
\caption{}\label{fig:bigon-S}
\end{figure}

Since $G_S$ has an $S$-cycle, $\widehat{T}$ is separating in $M(\beta)$, and so $\partial V_\beta$ is divided into two annuli $A_1$, $A_2$,
where $A_1$ meets $\partial D$.
Let $T_1=T\cup A_1$ and $T_2=T\cup A_2$.
Since $\partial D$ is non-separating on $T_1$, surgering $T_1$ along $D$ gives a torus disjoint from $K_\beta$.
On the other hand, surger $T_2$ along the face bounded by an $S$-cycle, which has distinct color from $D$.
This gives a torus disjoint from $K_\beta$.
Thus $M(\beta)$ is bounded by at most two tori, a contradiction.
Hence we can conclude that any component of $G_S^+$ containing a non-$S$-vertex is a cycle of bigon, possibly of length one.

Next, assume that $\Lambda$ has a cut vertex.
By Lemma \ref{lem:t2dblock}, there are only two possibilities for $\Lambda$ as shown in Figure \ref{fig:bigon-cycle2}.
However, we can still choose a disk face $D$ as in Figure \ref{fig:bigon-S}.
Thus a similar argument leads to a contradiction.
\end{proof}

\begin{figure}[tb]
\includegraphics*[scale=0.6]{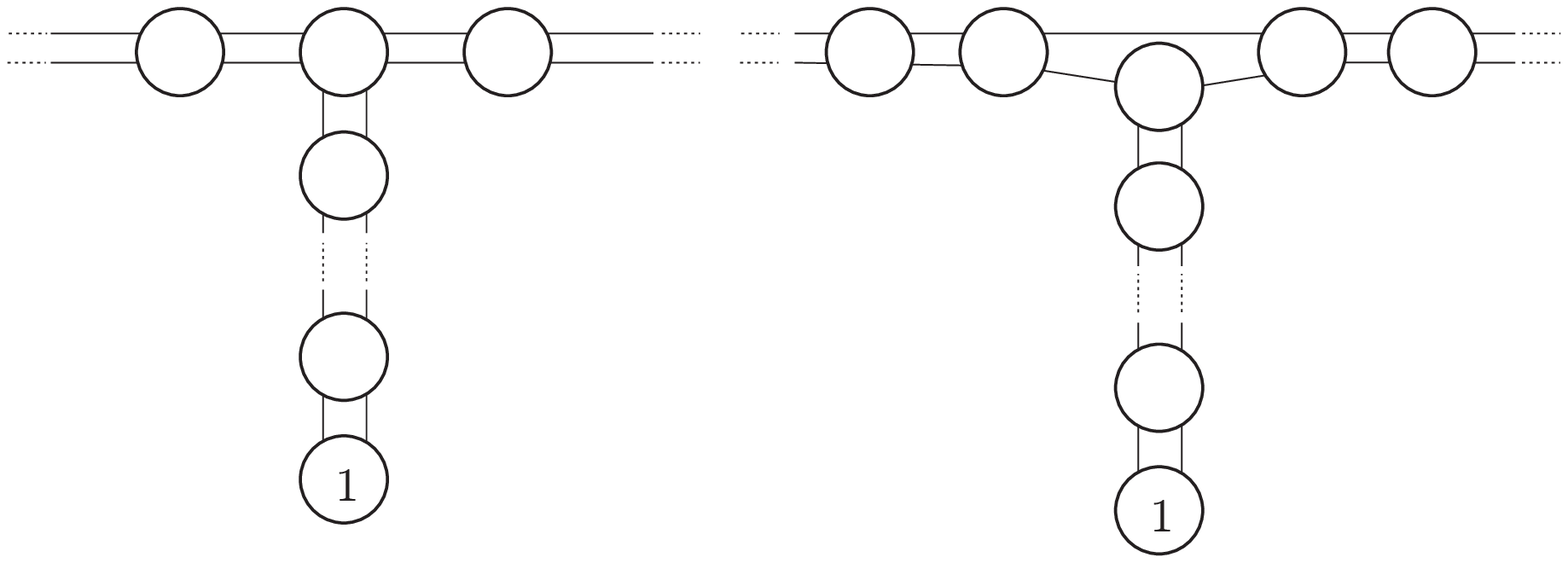}
\caption{}\label{fig:bigon-cycle2}
\end{figure}

\begin{proposition}
$t\ne 2$.
\end{proposition}

\begin{proof}
By Lemma \ref{lem:t2pre2}, any component of $G_S^+$ containing a non-$S$-vertex is a cycle of bigons.
All bigons have the same color by Lemma \ref{lem:key}, and hence 
any non-$S$-vertex is incident to exactly two adjacent negative edges.
This implies that each non-$S$-label appears once at each vertex of $G_T$ among positive edge endpoints.

Suppose that $G_T$ has no Scharlemann cycle.
Then every label appears once at each vertex among positive edge endpoints.
Also, two edges of the bigons belong to the same pair of families of mutually parallel negative edges in $G_T$ by \cite[Lemma 5.2]{GL3}.
(Otherwise, $M(\beta)$ would contain a Klein bottle meeting $K_\beta$ once.)
Hence $G_T$ has only two families of $s$ mutually parallel negative edges.
Thus $G_T\cong H'(s/2,s,s,0,0)$ or $H'(s/2,s,0,s,0)$.

If $G_T$ has a Scharlemann cycle, then each vertex of $G_T$ has at least $s+2$ positive edge endpoints,
and so just $s/2+1$ loops by Lemma \ref{lem:GT}(1),
two of which form an $S$-cycle.
Then we see that $G_T\cong H'(s/2+1,s,s-2,0,0)$ or $H'(s/2+1,s-2,0,s,0)$.

We consider these four cases.

Case (A): Assume $G_T\cong H'(s/2,s,s,0,0)$.
We can assume the labels in $G_T$ as in Figure \ref{fig:t2-1}(i).
Let $Q_1$ and $Q_2$ be the families of mutually parallel negative edges with $q_1(=s)$ and $q_2(=s)$ edges, respectively.
Let $\sigma$ be the associated permutation to $Q_1$ such that
an edge of $Q_1$ has label $x$ at $v_1$ and label $\sigma(x)$ at $v_2$.
Clearly, $Q_2$ also associates to the same permutation $\sigma$.
Since the edges of $Q_1$ and $Q_2$ form cycles of bigons in $G_S$, $\sigma^2$ is the identity.
Therefore $\sigma(x)=x$ or $\sigma(x)=x+s/2$.

\begin{figure}[tb]
\includegraphics*[scale=0.34]{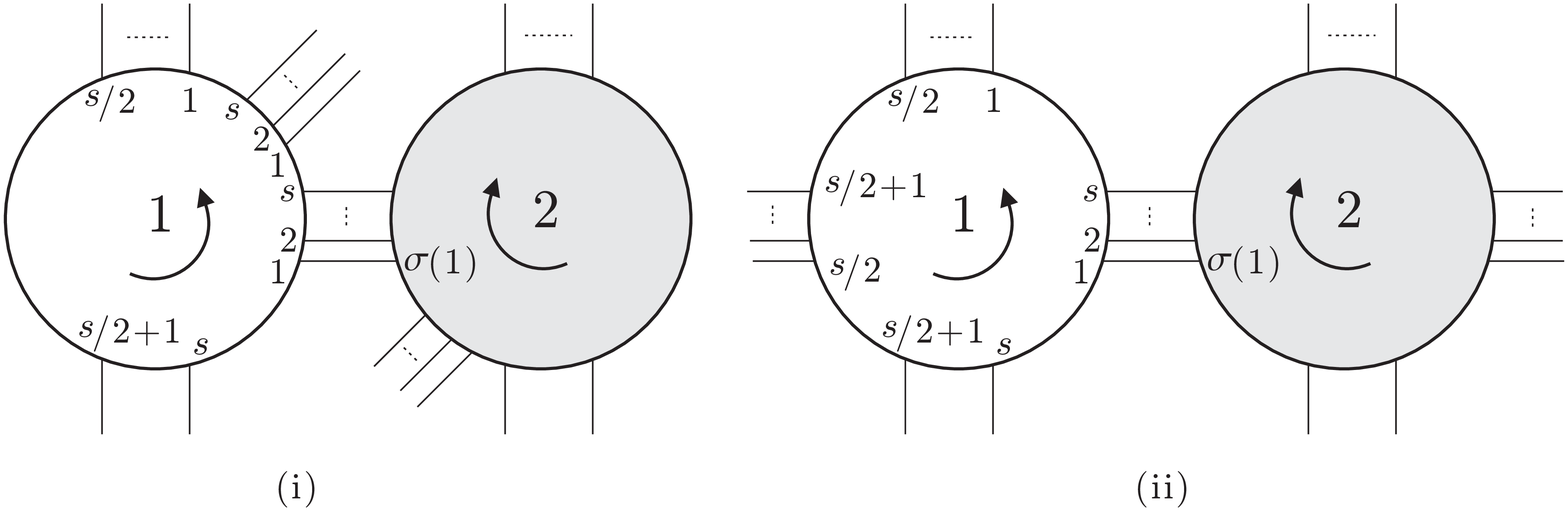}
\caption{}\label{fig:t2-1}
\end{figure}

Assume that $\sigma$ is the identity.
Then $G_S$ consists of $s/2$ copies of a graph isomorphic to $G(2,2,0)$ or $G(2,1,1)$.
Let $D$ be a $3$-gon in $G_T$.
Notice that $D$ is one-cornered.
Using $D$, one can see that $\widehat{S}$ is separating in $M(\alpha)$ and the side of $\widehat{S}$ containing $D$
is bounded by a torus.
Also, take a bigon $D'$ among the edges of $Q_1$, lying on the opposite side.
Moreover, we can choose $D'$ so that its edges bound an annulus in $\widehat{S}$ disjoint from the vertices of $G_S$,
as an innermost one.
Let $D'$ be bounded by an $x$-edge and a $(x+1)$-edge, and
let $A$ be the annulus in $\partial V_\alpha$ between $u_x$ and $u_{x+1}$.
Then surgering $(\widehat{S}-\mathrm{Int}\,(u_x\cup u_{x+1}))\cup A$ along $D'$ gives either an annulus, or
a disjoint union of an annulus and a torus.
In any case, the annulus component meets $K_\alpha$ fewer than $\widehat{S}$, and the torus component
is disjoint from $K_\alpha$.
Hence the annulus component is boundary parallel and the torus component is inessential.
Thus $M(\alpha)$ is bounded by at most two tori.

Next, assume that $\sigma(x)=x+s/2$.
Then we see that two $\{1,s\}$-loops in $G_T$ bound a bigon face $E$ in $G_S$.
But $\partial E$ runs like Figure \ref{fig:t2-11}(i), and so
$M(\beta)$ contains a Klein bottle meeting $K_\beta$ once, obtained from $E\cup H\cup A$ by shrinking $H$ radially into its core, where $H$ is the $1$-handle part of $V_\beta$
meeting $E$ and $A$ is the annular region on $\widehat{T}$ between the two $\{1,s\}$-loops.
This contradicts Theorem \ref{thm:noklein}.

\begin{figure}[tb]
\includegraphics*[scale=0.45]{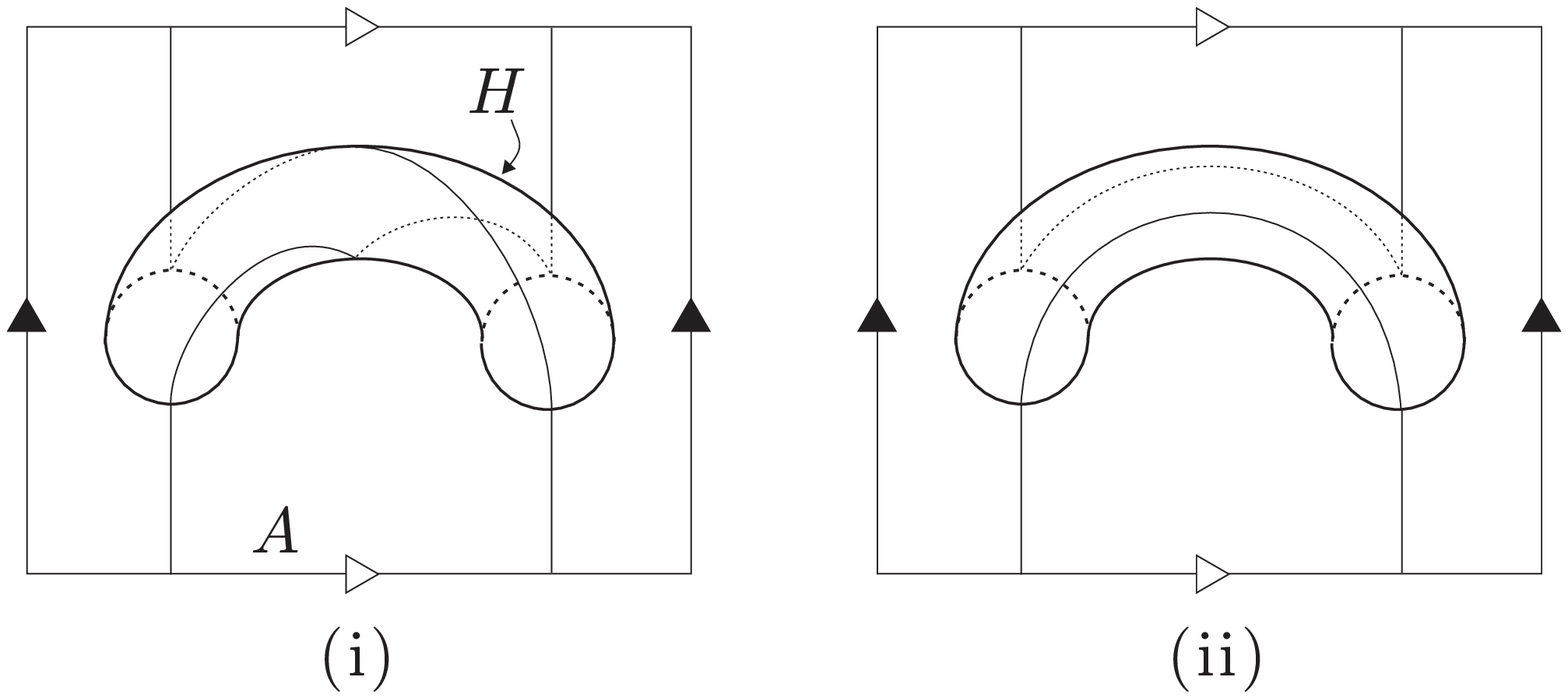}
\caption{}\label{fig:t2-11}
\end{figure}

Case (B): Assume $G_T\cong H'(s/2,s,0,s,0)$.
We can assume that the labels in $G_T$ are as in Figure \ref{fig:t2-1}(ii).
Similarly, we can see that two families $Q_1$ and $Q_2$ of mutually parallel negative edges associate to the same permutation $\sigma$,
and $\sigma^2$ is the identity.

If $\sigma$ is the identity, then take two $\{1,s\}$-loops in $G_T$.
They bound a bigon $E$ in $G_S$, and $\partial E$ runs like Figure \ref{fig:t2-11}(ii).
But consider any $S$-cycle in $G_S$.
It has one edge in each of $Q_1$ and $Q_2$, but we cannot connect them on $\partial V_\beta$.

When $\sigma(x)=x+s/2$, the same argument as in Case (A) gives a contradiction.

Case (C): Assume $G_T\cong H'(s/2+1,s,s-2,0,0)$.
The labels in $G_T$ can be assumed as in Figure \ref{fig:t2-2}(i).
But this implies that the component of $G_S^+$ containing $u_3$ is not a cycle of bigons, contradicting Lemma \ref{lem:t2pre2}.

\begin{figure}[tb]
\includegraphics*[scale=0.34]{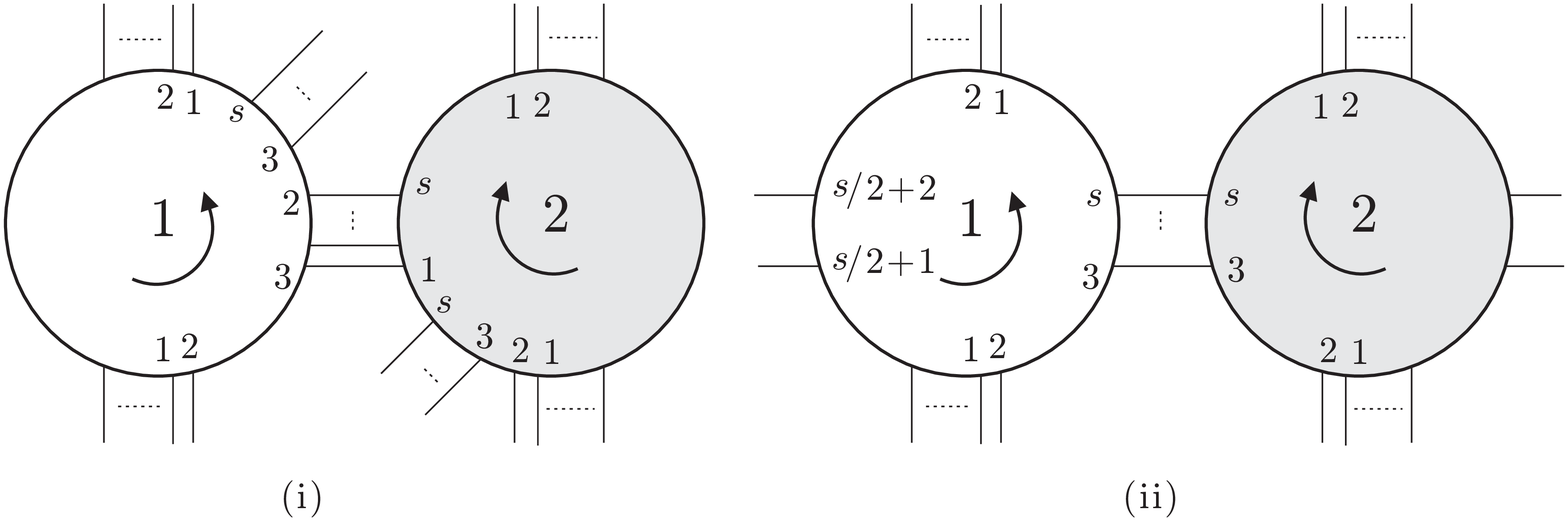}
\caption{}\label{fig:t2-2}
\end{figure}

Case (D): Assume $G_T\cong H'(s/2+1,s-2,0,s,0)$.
Then the labels in $G_T$ can be assumed as in Figure \ref{fig:t2-2}(ii).
Two $\{3,s\}$-loops in $G_T$ bound a bigon in $G_S$.
Then the same argument as in Case (B) leads to a contradiction.
\end{proof}

\subsection{Cases $t\ge 5$ or $t=3$}\label{subsec:nonS}

In this subsection, we eliminate two cases where $t\ge 5$ and $t=3$.
When $t=3$, $G_S$ has no Scharlemann cycle by Lemma \ref{lem:common}(4), and so no $S$-label.
If $t\ge 5$, then $G_S$ has a non-$S$-label by Lemma \ref{lem:GS}(2).

\begin{lemma}\label{lem:dblock}
$G_S^+$ has no disk block.
\end{lemma}

\begin{proof}
Let $B$ be a disk block of $G_S^+$.
It has at most one cut vertex of $G_S^+$ and at most one $S$-vertex among boundary vertices.
Let $V_i$, $V_b$, $V_c$, $V_s$ be the number of interior, boundary, cut and $S$-vertices, respectively.
Here a cut vertex means a cut vertex of $G_S^+$.
Then $V_c, V_s\le 1$.
Possibly, an $S$-vertex is a cut vertex.
In this case, we set $V_s=0$ and $V_c=1$.

Let $x$ be a non-$S$-label.
Any interior vertex is incident to three positive $x$-edges, any boundary vertex, except
a cut vertex and an $S$-vertex, is incident to at least two such edges by Lemma \ref{lem:2t},
and an $S$-vertex is incident to at least one such edge by Lemma \ref{lem:t}.
Consider the subgraph $B^x$ of $B$ consisting of all vertices and all $x$-edges of $B$.
We remark that $B^x$ may be disconnected, and may have many cut vertices.
Let $V$, $E$, $F$ be the number of vertices, edges, disk faces of $B^x$, respectively, as a graph in a disk.
Then $V=V_i+V_b$ and $F\ge 1-V+E$.
By counting $x$-edges, we have
\begin{equation}
E\ge 3V_i+2(V_b-V_c-V_s)+V_s=3V-V_b-2V_c-V_s. \label{eqn:x-edge}
\end{equation}
Since each disk face of $B^x$ has at least $4$ sides by Lemma \ref{lem:key2},
\begin{equation}
2E\ge 4F+V_b'\ge 4(1-V+E)+V_b, \label{eqn:face}
\end{equation}
where $V_b'$ is the number of boundary vertices of $B^x$.
(Notice $V_b'\ge V_b$.)
These give
$3V-V_b-2V_c-V_s\le 2V-2-V_b/2$.
Equivalently, $V-V_b/2+2\le 2V_c+V_s$.
Hence $V_c=V_s=1$ and $V=V_b=2$.
This implies that $B$ is a family of at least $t$ parallel positive edges joining two vertices, which contradicts Lemma \ref{lem:parallel-max}(2).
\end{proof}

\begin{lemma}\label{lem:a-ex}
Any component of $G_S^+$ has an annulus support, and is extremal.
\end{lemma}

\begin{proof}
If $G_S^+$ has a component with a disk support,
then there is an extremal one $\Lambda$ with a disk support.
Hence $\Lambda$ contains a disk block, contradicting Lemma \ref{lem:dblock}.
\end{proof}

\begin{proposition}\label{prop:nonS}
$t=4$.
\end{proposition}

\begin{proof}
Choose an outermost component $\Lambda$ of $G_S^+$.
There is an annulus $A$ in $\widehat{S}$ such that $\Lambda\subset \mathrm{Int}\,A$,
$A\cap G_S^+=\Lambda$ and $A$ contains one component of $\partial \widehat{S}$.
After capping off that component of $\partial \widehat{S}$ with a disk,
we regard $\Lambda$ as lying in a disk.
In this view point, we consider its interior and boundary vertices.
Let $V_i$, $V_b$, $V_s$ be the number of interior, boundary, and $S$-vertices of $\Lambda$, respectively.
Remark that $\Lambda$ has a disk face $f$ containing the disk capped off in its interior, where $f$ may be a monogon.
Also, $\Lambda$ may have an $S$-vertex, and
a cut vertex (of $\Lambda$) among boundary vertices.
But any boundary vertex is good by Lemma \ref{lem:dblock}.

Let $x$ be a non-$S$-label.
Consider a subgraph $\Lambda^x$ of $\Lambda$ consisting of all vertices and all $x$-edges of $\Lambda$, as a graph in a disk.
We remark that $\Lambda^x$ may be disconnected.
Let $V$, $E$, $F$ be the number of vertices, edges, disk faces of $\Lambda^x$.
Then $F\ge 1-V+E$ and $V=V_i+V_b$.
Each interior vertex of $\Lambda$ is incident to three positive $x$-edges,
each boundary vertex is incident to at least two such edges, and an $S$-vertex is incident to at least one such edge.
Hence we have
\begin{equation}
E\ge 3V_i+2(V_b-V_s)+V_s=3V-V_b-V_s. \label{eqn:a-x-edge}
\end{equation}
Also, since each disk face of $\Lambda^x$, possibly except one, has at least $4$ sides,
\begin{equation}
2E\ge 4(F-1)+1+V_b'\ge 4(E-V)+1+V_b, \label{eqn:a-face}
\end{equation}
where $V_b'$ is the number of boundary vertices of $\Lambda^x$ itself.
These give
$3V-V_b-V_s\le E\le 2V-V_b/2-1/2$,
equivalently $V_i+V_b/2\le V_s-1/2$.
Then $V_s=1$, $V_i=0$ and $V_b=1$.
This means that $\Lambda$ is an $S$-vertex with at least $t/2$ parallel loops.

Choose another outermost component $\Lambda'$ of $G_S^+$ near the other component of $\partial \widehat{S}$.
The same argument shows that $\Lambda'$ consists of an $S$-vertex  and at least $t/2$ parallel loops.
Since two $S$-vertices are connected with the edges of Scharlemann cycles,
$G_S^+$ cannot have other components than $\Lambda$ and $\Lambda'$.
But this means $s=2$, a contradiction.
\end{proof}

\subsection{Case $t=4$}\label{subsec:t4}

Again, we can show that $G_S^+$ has no disk block as in Lemma \ref{lem:dblock}, but it needs another argument.

\begin{lemma}\label{lem:dblock2}
$G_S^+$ has no disk block.
\end{lemma}

\begin{proof}
Let $B$ be a disk block.
We use the same notation as in the proof of Lemma \ref{lem:dblock}.
By Lemma \ref{lem:GS}(1), we can choose a label $x$ which is not a label of an $S$-cycle.
Then \eqref{eqn:x-edge} holds.
Since each disk face of $B^x$ has at least three sides, \eqref{eqn:face} changes to
$2E\ge 3F+V_b'\ge 3(1-V+E)+V_b$.
These give
$3V-V_b-2V_c-V_s\le E\le 3V-3-V_b$,
equivalently, $2V_c+V_s\ge 3$.
Hence $V_c=V_s=1$, and all inequalities above are equalities.
So, $E=3V-3-V_b$, $F=2V-V_b-2=2V_i+V_b-2\ge 0$ and each disk face of $B^x$ is $3$-sided.
If $F=0$, then $V=V_b=2$ and hence $B$ is a family of at least $t$ mutually parallel edges, contradicting Lemma \ref{lem:parallel-max}(2).
Thus $F>0$.

We may assume that $x=4$ without loss of generality.
Figure \ref{fig:3side} lists all possible $3$-sided faces of $B^x$,
where all edges of $G_S$ are indicated.

\begin{figure}[tb]
\includegraphics*[scale=0.65]{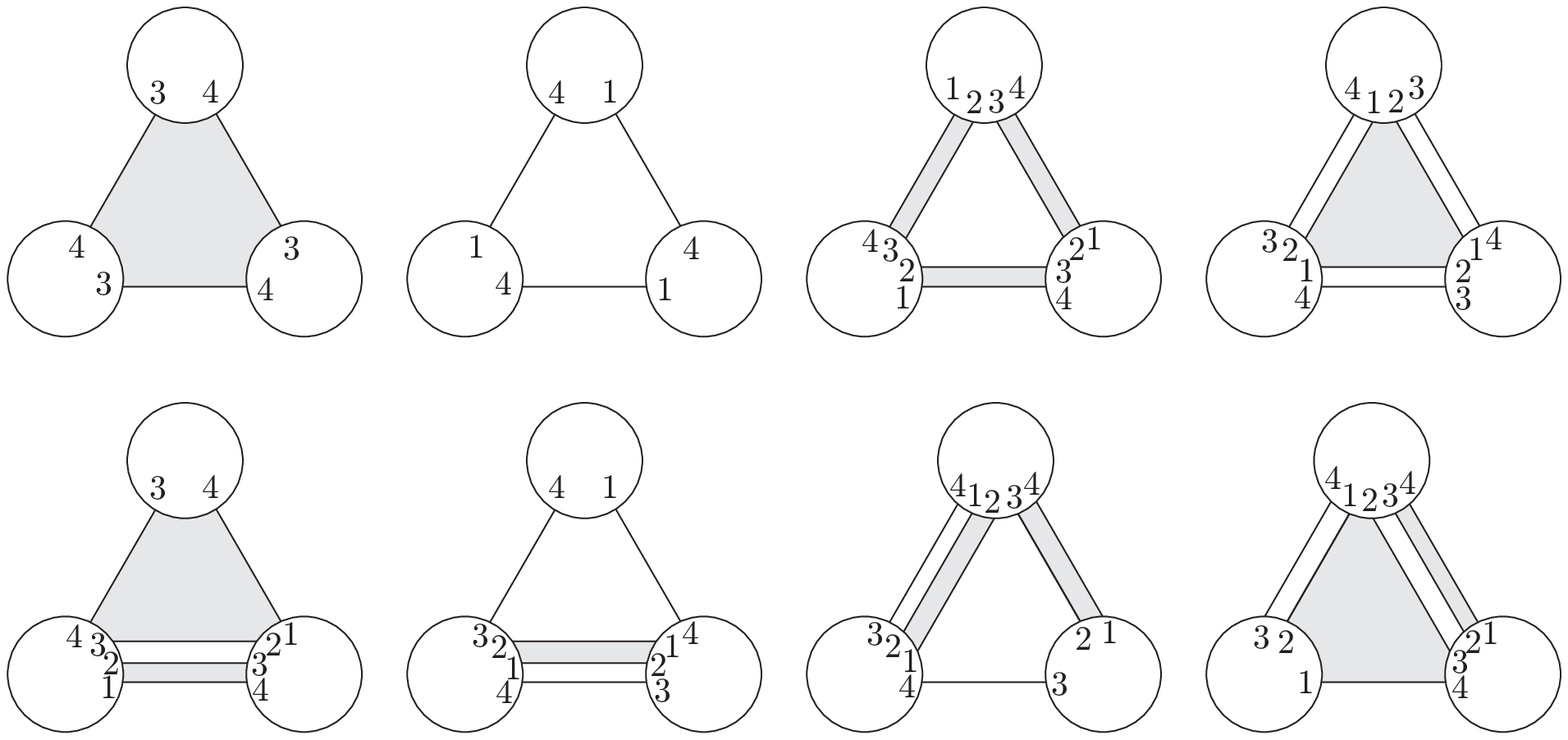}
\caption{}\label{fig:3side}
\end{figure}

Extended Scharlemann cycles are impossible.
The last four configurations can be eliminated in the same way.
For example, it contains a black $S$-cycle and two two-cornered white faces, a bigon and a $3$-gon adjacent to the $S$-cycle.
These white faces are homologically independent.
Hence $M(\beta)$ is bounded by at most two tori.
Thus only the first and second configurations are possible, and they cannot occur simultaneously by Lemma \ref{lem:key}.
Hence we may assume that all faces of $B^x$ are bounded by black Scharlemann cycles with label pair $\{3,4\}$.
Of course, this is impossible if $F>1$.
But if $F=1$, then $V_i=V_b=1$ or $V=V_b=3$.
In the former, $B^x$ has a vertex of degree one, which contradicts Lemma \ref{lem:parallel-max}(2).
In the latter, $B^x$ is a cycle of length three, and so the vertex other than the cut vertex and the $S$-vertex is incident to
at least $t$ parallel positive edges in $B$, contradicting Lemma \ref{lem:parallel-max}(2) again.
\end{proof}

Hence Lemma \ref{lem:a-ex} holds again.

\begin{proposition}
$t\ne 4$.
\end{proposition}

\begin{proof}
Assume $t=4$.  We use the same notation as in the proof of Proposition \ref{prop:nonS}.
Let $x$ be a label of $G_S$ which is not a label of an $S$-cycle.
Then we have \eqref{eqn:a-x-edge}.
Since each disk face of $\Lambda^x$, possibly except one, has at least three sides, \eqref{eqn:a-face} changes to
$2E\ge 3(F-1)+1+V_b'\ge 3(E-V)+1+V_b$,
where $V_b'$ denotes the number of boundary vertices of $\Lambda^x$.
These give
$3V-V_b-V_s\le E\le 3V-V_b-1$.
Hence $V_s=1$ and all inequalities above are equalities, and then
$E=3V-V_b-1$, $F=2V-V_b=2V_i+V_b$, and each disk face of $\Lambda^x$ is $3$-sided.

If $F=2$, then $V=V_b=2$ and $E=3$.
Then $\Lambda$ has two vertices, one of which is a pinched vertex and the other is an $S$-vertex.
By examining the labels around the vertices, we can see that $\Lambda$ contains two $S$-cycles with disjoint label pairs.
This contradicts Lemma \ref{lem:parallel-max}(1).
If $F>2$, then the same argument as in the proof of Lemma \ref{lem:dblock2} is applicable.
Thus $F=1$.
Then $V=V_b=1$ and so $\Lambda$ consists of an $S$-vertex and parallel loops.

Similarly, another outermost component of $G_S^+$ near the other component of $\partial \widehat{S}$ consists
of an $S$-vertex with parallel loops.
Then $s=2$ as in the proof of Proposition \ref{prop:nonS}, a contradiction.
\end{proof}

\section{Klein bottle}\label{sec:klein}

In the rest of paper, we prove Theorem \ref{thm:noklein}.
Suppose that $M(\beta)$ contains a Klein bottle $\widehat{P}$ which meets $K_\beta$ in $p\ (\le t/2)$ points,
and that $p$ is minimal among all Klein bottles in $M(\beta)$.
Then $\widehat{P}$ meets $V_\beta$ in a disjoint union of meridian disks $w_1,w_2,\dots,w_p$ numbered successively along
$V_\beta$.
Let $P=\widehat{P}\cap M$, and let $N$ be a thin neighborhood of $\widehat{P}$.

\begin{lemma}
$P$ is incompressible and boundary incompressible.
\end{lemma}

\begin{proof}
See \cite[Lemma 2.1]{LOT}.
\end{proof}

Thus we can assume that no circle component of $S\cap P$ bounds a disk in $S$ or $P$.
From the arc components of $S\cap P$, we have a graph pair in the usual way.
By abuse of notation, we denote the pair by $\{G_S,G_P\}$ in the rest of paper.
Since $P$ is non-orientable, we cannot give a sign to a vertex of $G_P$.
However, there is a way to give a sign to an edge of $G_P$ (see \cite{JLOT}).
Then the parity rule survives without any change.
Remark that a positive edge of $G_S$ can be an level edge.
It corresponds to an orientation-reversing loop on $\widehat{P}$.
Also, there are no two edges which are parallel in both graphs \cite[Lemma 2.1]{Go2}.

If $p>2$, 
a triple $\{e_1,e_2,e_3\}$ of mutually parallel positive edges in $G_S$ is called
a \textit{generalized $S$-cycle\/} if $e_2$ is a level edge with label $i$,
and $e_1$ and $e_3$ have label pair $\{i-1,i+1\}$ at their endpoints.

\begin{lemma}\label{lem:parallel-kb}
If $p\ge 2$, then $G_S$ satisfies the following.
\begin{itemize}
\item[(1)] There is no Scharlemann cycle.
\item[(2)] If $p\ge 3$, then there is no generalized $S$-cycle.
\item[(3)] At most two labels can be labels of positive level edges.
\item[(4)] Any family of parallel positive edges contains at most $p/2+1$ edges.
\item[(5)] Any family of parallel negative edges contains at most $p$ edges.

\end{itemize}
\end{lemma}

\begin{proof}
(1) See \cite[Lemma 3.2]{T2}.  (The argument works for a Scharlemann cycle with any length.)
(2) is \cite[Lemma 3.3]{T2}.
(3) follows from the facts that a positive level edge in $G_S$ corresponds to an orientation-reversing loop in $\widehat{P}$ and that
a Klein bottle contains at most two disjoint M\"{o}bius bands.

(4) Let $Q$ be the family of mutually parallel positive edges in $G_S$.
Let $|Q|$ denote the number of edges in $Q$.  Suppose $|Q|>p/2+1$.

Assume $p=2$.
If an edge in $Q$ is level, then all edges are level.
Since any two level edges with the same label are parallel in $G_P$,
there would be two edges which are parallel in both graphs, a contradiction.
If no edge in $Q$ is level, then $Q$ contains an $S$-cycle, contradicting (1).

Assume $p>2$.
Then $Q$ would contain an $S$-cycle or a generalized $S$-cycle, a contradiction. 


(5) Let $e_1,e_2,\dots,e_p,e_1'$ be mutually parallel negative edges in $G_S$, numbered successively.
We may assume that $e_i$ has label $i$ at one vertex for $i=1,2,\dots,p$, so $e_1'$ has label $1$ at the same vertex.
If $e_i$ has label $\sigma(i)$ at the other end, we have the associated permutation $\sigma$.
According to the orbits of $\sigma$, the edges $e_i$ form essential orientation-preserving cycles on $\widehat{P}$ by \cite[Lemma 2.3]{Go2}.
Let $L$ be the cycle through vertex $w_1$.
Then $e_1'$ is not parallel to $e_1$.
However then a new cycle $(L-e_1)\cup e_1'$ is inessential on $\widehat{P}$, a contradiction.
(This is essentially the same as the proof of \cite[Lemma 4.2]{Go2}.)
\end{proof}

\begin{lemma}\label{lem:key3}
Let $p\ge 3$.
If $x$ is not a label of a positive level edge in $G_S$,
then any $x$-face in $G_S$ has at least $4$ sides.
\end{lemma}

\begin{proof}
First, there is no two-sided $x$-face, since it contains an $S$-cycle or a generalized $S$-cycle.
Let $D$ be a $3$-sided $x$-face, and let $\Gamma=G_S\cap D$.
If $\Gamma$ does not contain a level edge, then there is a Scharlemann cycle by \cite{HM}, contradicting Lemma \ref{lem:parallel-kb}(1).
Hence $\Gamma$ contains a level edge.
Notice that the faces of $\Gamma$ consist of a single $3$-gon $f$ and bigons.
Since $\Gamma$ cannot contain a generalized $S$-cycle, any level edge appears in the $3$-gon $f$.
There are two cases.

(1) Only one label is a label of positive level edges in $\Gamma$.

Then, in fact, $\Gamma$ contains only one level edge $e$.
We may assume that it has label $1$.
Clearly, the bigon $g$ adjacent to $e$ has two corners $(1,2)$ and $(p,1)$.
Moreover, the $3$-gon $f$ is also two-corned.
That is, it has only $(1,2)$-corner and $(p,1)$-corner \cite[Claim 3.7]{IT} (or see \cite{LOT}).

Let $H$ be the part of $V_\beta$ between $w_p$ and $w_1$, containing $w_2$.
Let $X=N(\widehat{P}\cup H\cup f\cup g)$.
Then $\partial X$ is a torus intersecting $K_\beta$ fewer than $t$ times.
Hence it is boundary parallel in $M(\beta)$ or compressible.
Thus $M(\beta)$ is bounded by at most one torus, a contradiction.

(2) Two labels are labels of positive level edges in $\Gamma$.

We may assume that the $3$-gon $f$ contains a level edge $e_1$ with label $1$ and a level edge $e_2$ with label $2$.
Let $g_i$ be the bigon adjacent to $f$, sharing $e_i$ for $i=1,2$.
Let $H$ be the part of $V_\beta$ between $w_p$ and $w_3$, containing $w_1$.
Construct $N(\widehat{P} \cup H \cup f\cup g_1\cup g_2)$ as above.
Then a similar argument to (1) implies a contradiction.
\end{proof}

\begin{lemma}\label{lem:s1}
$s\ne 1$.
\end{lemma}

\begin{proof}
Assume $s=1$.
Notice that $p$ is even, since the vertex of $G_S$ has degree $3p$.
There are $3p/2$ parallel loops in $G_S$, but this contradicts Lemma \ref{lem:parallel-kb}(4), because $3p/2>p/2+1$.
\end{proof}




\begin{lemma}\label{lem:p1}
$p\ne 1$.
\end{lemma}

\begin{proof}
Assume $p=1$.  
By an Euler characteristic calculation, $G_S$ has a disk face $D$.
Let $X=N\cup V_\beta$.
Then $\partial X$ is a genus two closed surface disjoint from $K_\beta$.
Let $D'=D-\mathrm{Int}\, X$.
Surger $\partial X$ along $D'$.
The resulting surface is either a torus or a disjoint union of two tori, according as $\partial D'$ is non-separating or separating on $\partial X$.
Thus $M(\beta)$ is bounded by at most two tori, a contradiction.
\end{proof}

\begin{lemma}
$s\ge 3$.
\end{lemma}

\begin{proof}
By Lemmas \ref{lem:s1} and \ref{lem:p1}, $s\ge 2$ and $p\ge 2$.
Suppose $s=2$.
Then we can use the same notation $G_S\cong G(p_0,p_1,p_2)$ as in Section \ref{sec:two}.

First assume $p=2$.
Since $G_S$ cannot contain an $S$-cycle, $p_0\le 1$.
By Lemma \ref{lem:parallel-kb}(4) and (5), $p_i\le 2$ for $i=1,2$.
Thus $G_S\cong G(1,2,2)$, and there are two bigons and two $3$-gons.
Take a bigon $D_1$ and a $3$-gon $D_2$.
Let $X=N\cup V_\beta$, and $D_i'=D_i-\mathrm{Int}\,X$.
Then $\partial X$ is a genus three closed surface, on which $\partial D_1'$ and $\partial D_2'$ are
homologically independent.
Thus $M(\beta)$ is bounded by at most one torus, a contradiction.

Next, assume $p\ge 3$.
By Lemma \ref{lem:parallel-kb}(4), $p_0\le p/2+1$, but
if the equality holds, there is an $S$-cycle.
Hence $p_0\le (p+1)/2$.

If the two vertices of $G_S$ are parallel, then $p_i\le p/2+1$ for $i=1,2$.
So $3p\le 2\cdot (p+1)/2+2(p/2+1)=2p+3$, giving $p\le 3$.
When $p=3$, we have $p_i\le 2$, giving $3p\le 2\cdot (p+1)/2+2\cdot(p+1)/2=2p+2$.
This is a contradiction.

Therefore the two vertices of $G_S$ are antiparallel.
By Lemma \ref{lem:parallel-kb}(5), $p_i\le p$ for $i=1,2$.
So $3p=2p_0+p_1+p_2\le 2p_0+2p$, giving $p_0\ge p/2$.
Hence $p_0=p/2$ if $p$ is even, and $p_0=(p+1)/2$ if $p$ is odd.
This implies that $G_S\cong G(p/2,p,p)$ if $p$ is even, and $G_S\cong G((p+1)/2,p,p-1)$ if $p$ is odd.
For both cases, the same argument as the case $G_S\cong G(t/2,t,t)$ in the proof of Lemma \ref{lem:t2} works.
\end{proof}

\begin{lemma}\label{lem:2p-p}
Any vertex of $G_S$, except $S$-vertices, has at least $2p$ positive edge endpoints.
An $S$-vertex, if it exists, has at least $p$ positive edge endpoints.
\end{lemma}

\begin{proof}
Lemma \ref{lem:xface} holds again.
Hence the proofs of Lemmas \ref{lem:2t} and \ref{lem:t} work.
\end{proof}


\begin{proposition}\label{prop:pge3}
$p\ge 3$ is impossible.
\end{proposition}

\begin{proof}
Using Lemma \ref{lem:key3}, instead of Lemma \ref{lem:key2}, the proof of Lemma \ref{lem:dblock}
works.
Hence $G_S^+$ does not contain a disk block.
Then the proofs of Lemma \ref{lem:a-ex} and Proposition \ref{prop:nonS} are applicable.
\end{proof}

\section{A special case: $p=2$}\label{sec:p=2}

Finally, we eliminate the situation where
$s\ge 3$ and $p=2$.
Recall that any vertex of $G_S$, except $S$-vertices, has at least four positive edge endpoints, and 
that any $S$-vertex, if it exists, has at least two positive edge endpoints by Lemma \ref{lem:2p-p}.


Let $W=\mathrm{cl}(M(\beta)-N)$.
We say that $N$ is a black region, and $W$ is a white region.
Let $T=\partial N-\mathrm{Int}\, V_\beta$.
As usual, $S$ and $T$ give a labelled graph pair $\{G_S',G_T\}$.
In fact, $G_T$ is a double cover of $G_P$.
The disk faces of $G_S'$ are divided into black and white faces as usual.
Thus any black bigon of $G_S'$ corresponds to an edge of $G_S$.

Consider a genus three surface $R=\partial(N\cup V_\beta)$, which is disjoint from $K_\beta$.

\begin{lemma}\label{lem:white}
For any two white disk faces of $G_S'^+$, their boundaries are parallel in $R$.
In particular, all white disk faces of $G_S^+$ have the same number of sides, and $G_S^+$ cannot contain
two adjacent $3$-gons.
\end{lemma}

\begin{proof}
Suppose that $G_S'^+$ contains two white disk faces whose boundaries are not parallel in $R$.
Surgering $R$ along them gives a torus or a disjoint union of two tori.
Since the surface is disjoint from $K_\beta$, $M(\beta)$ is bounded by at most two tori, a contradiction.

Let $f$ and $g$ be adjacent $3$-gons in $G_S^+$.
Consider two white faces $f'$ and $g'$ of $G_S'^+$ corresponding to $f$ and $g$, respectively.
Then $\partial f'$ and $\partial g'$ are not parallel on $R$.
\end{proof}

\begin{lemma}\label{lem:twobigon}
At any vertex of $G_S$, there are no consecutive pairs of parallel positive edges.
\end{lemma}

\begin{proof}
Otherwise, there are two consecutive bigons.
However, it is easy to see that the corresponding white bigons have non-parallel boundaries on $R$.
This contradicts Lemma \ref{lem:white}.
(See also \cite[Lemma 6.3]{JLOT}.)
\end{proof}

We divide the argument into two cases.

Case (A): $G_S^+$ contains a bigon.
Then all disk faces of $G_S^+$ are bigons by Lemma \ref{lem:white}.

\begin{lemma}\label{lem:form-dblock}
$G_S^+$ has at most two disk blocks.
Any disk block consists of two vertices, one of which is an $S$-vertex, and a pair of parallel edges.
\end{lemma}

\begin{proof}
Let $B$ be a disk block.
Since all faces of $B$ are disks, they are bigons.
Thus $B$ has only two vertices.
By Lemma \ref{lem:twobigon}, one vertex is an $S$-vertex.
Also, the other is a cut vertex of $G_S^+$.

Since $G_S$ has at most two $S$-vertices, there are at most two disk blocks in $G_S^+$.
\end{proof}

\begin{lemma}\label{lem:p2asupp}
Any component of $G_S^+$ has an annulus support, and is extremal.
\end{lemma}

\begin{proof}
If there is a component with a disk support, then there is an extremal one $\Lambda$ with a disk support.
Notice that all faces of $\Lambda$ are disk, and hence bigons. 
Thus there are two consecutive bigons at a non-$S$-vertex, which is not a cut vertex of $\Lambda$.
This contradicts Lemma \ref{lem:twobigon}.
\end{proof}

\begin{lemma}\label{lem:out}
Let $\Lambda$ be an outermost component of $G_S^+$.
Then $\Lambda$ consists of two vertices together with a loop at one vertex and a pair of parallel level edges connecting two vertices.
Moreover, one vertex is an $S$-vertex.
\end{lemma}

\begin{proof}
By Lemma \ref{lem:twobigon}, $\Lambda$ has no interior vertex.
If $\Lambda$ has no disk block, then it is a cycle of bigons, contradicting Lemma \ref{lem:twobigon}.
(If the cycle is length one, then there is an $S$-cycle.)
Also, any boundary vertex is incident to a disk block.
Since there is only one disk block incident to $\Lambda$,
we have the conclusion.
\end{proof}

\begin{lemma}\label{lem:caseA}
Case \textup{(A)} is impossible.
\end{lemma}

\begin{proof}
By Lemmas \ref{lem:p2asupp} and \ref{lem:out},
$G_S^+$ has two components $\Lambda_1$ and $\Lambda_2$, each of which satisfies the conclusion of Lemma \ref{lem:out}.
We may assume that $\Lambda_i$ contains an $S$-vertex $u_i$ for $i=1,2$.
Notice that $u_1$ and $u_2$ are joined by the edges of a Scharlemann cycle in $G_T$, which do not lie on a disk in $\widehat{S}$
by Lemma \ref{lem:common}(3).
Hence $G_S^+$ consists of $\Lambda_1$ and $\Lambda_2$, so $s=4$.
Since $u_1$ is incident to four negative edges, $G_P$ contains at least two $1$-faces by Euler characteristic calculation.
Each $1$-face contains a Scharlemann cycle.
Thus $G_P$ has at least two Scharlemann cycles, so the $4$ negative edges at $u_1$ are the edges of Scharlemann cycles in $G_P$.
This is similar for $u_2$.
Then the non-$S$-vertex of $\Lambda_1$ cannot be incident to a negative edge, a contradiction.
\end{proof}

Case (B): Any disk face of $G_S^+$ has at least three sides.

\begin{lemma}\label{lem:p2dblock}
$G_S^+$ has no disk block.
\end{lemma}

\begin{proof}
Let $B$ be a disk block.
It has at most one cut vertex and at most one $S$-vertex among boundary vertices.
Let $V$, $E$, $F$ be the number of vertices, edges, faces of $B$, respectively.
Let $V_i$, $V_b$, $V_c$, $V_s$ be the number of interior, boundary, cut, and $S$-vertices of $B$, respectively.
Then $V=V_i+V_b$ and $V_c, V_s\le 1$.
(If an $S$-vertex is a cut vertex, then set $V_c=1$ and $V_s=0$.)

Any interior vertex has degree $6$, any boundary vertex, except a cut vertex and an $S$-vertex, has degree at least $4$,
and a cut vertex or an $S$-vertex has degree at least two. 
By counting degree,
$2E\ge 6V_i+4(V_b-V_c-V_s)+2V_c+2V_s=6V-2V_b-2V_c-2V_s$.
Since each face of $B$ has at least three sides,
$2E\ge 3F+V_b=3(1-V+E)+V_b$.
Then $3V-V_b-V_c-V_s\le 3V-3-V_b$, and hence $V_c+V_s\ge 3$, a contradiction.
\end{proof}

\begin{lemma}
Case \textup{(B)} is impossible.
\end{lemma}

\begin{proof}
By Lemma \ref{lem:p2dblock},
any component of $G_S^+$ has an annulus support, and is extremal.
Let $\Lambda$ be an outermost component.
After capping off the component of $\partial \widehat{S}$ near $\Lambda$,
we regard $\Lambda$  as lying in a disk.
In this view points, we consider its interior vertices and boundary vertices.
Let $V$, $E$, $F$ be the number of vertices, edges, and disk faces of $\Lambda$, respectively.
Let $V_i$, $V_b$, $V_s$ be the number of interior, boundary, and $S$-vertices of $\Lambda$.
Here $\Lambda$ may have a monogon, which includes the disk capped off.
As before, $2E\ge 6V_i+4(V_b-V_s)+2V_s=6V-2V_b-2V_s$.
Since each disk face of $\Lambda$, except at most one, has at least three sides,
$2E\ge 3(F-1)+1+V_b=3E-3V+1+V_b$.
Then 
$3V-V_b-V_s\le 3V-V_b-1$, equivalently, $V_s\ge 1$.
Thus $V_s=1$ and all inequalities above are equalities.
So, each disk face of $\Lambda$, except one monogon, has three sides.
Since $\Lambda$ has an $S$-vertex, $\widehat{S}$ is separating in $M(\alpha)$ and $G_S^+$ has exactly two components, $\Lambda$ and $\Lambda'$,
where $\Lambda'$ is another outermost component.

If $F=1-V+E=2V-V_b>2$, then 
$\Lambda$ contains two adjacent $3$-gons, contradicting Lemma \ref{lem:white}.
If $F=1$, then $V=V_b=V_s=1$ and $E=1$.
Hence $\Lambda$ is an $S$-vertex with a loop.
Similarly, $\Lambda'$ has the same form.
But this means $s=2$, a contradiction. 
If $F=2$, then $V=V_b=2$ and $E=3$.
Then $\Lambda$ consists of one pinched vertex and one bivalent $S$-vertex.
Again, $\Lambda'$ has the same form.
Since $u_1$ is incident to four negative edges, $G_P$ contains at least two Scharlemann cycle as in the proof of Lemma \ref{lem:caseA}.
Then any pinched vertex cannot be incident to a negative edge, a contradiction.
\end{proof}

\bibliographystyle{amsplain}

\end{document}